\documentclass[11pt,twoside]{article}

\title{Vector bundles near negative curves:\\moduli and local Euler characteristic}
\author{E.\ Ballico \and E.\ Gasparim \and T.\ K\"{o}ppe}
\date{}

\usepackage[pdftitle={Vector bundles near negative curves: moduli and local Euler characteristic},
            pdfauthor={Edoardo Ballico, Elizabeth Gasparim, Thomas Köppe},
            pdfsubject={Mathematics, Algebraic Geometry}, pdfproducer={}, pdfcreator={},
            pdfpagelayout=TwoColumnRight, colorlinks=true, linkcolor=black, urlcolor=black, citecolor=black]{hyperref}
\usepackage{amsmath,amssymb,amsthm,calc,fancyhdr,mathrsfs,cancel}
\usepackage[outer=3cm,inner=3cm]{geometry}

\numberwithin{equation}{section}
\newtheorem{theorem}{Theorem}[section]
\newtheorem{lemma}[theorem]{Lemma}
\newtheorem{corollary}[theorem]{Corollary}
\newtheorem{proposition}[theorem]{Proposition}
\newtheorem*{theorem*}{Theorem}
\newtheorem*{thm_genopset}{Theorem \ref{thm.genopset}}
\newtheorem*{thm_hdcomps}{Theorem \ref{thm.hdcomps}}
\newtheorem*{thm_admissible}{Theorem \ref{thm.admissible}}
\newtheorem*{thm_height_nonsplit}{Theorem \ref{thm.height_nonsplit}}
\newtheorem*{thm_bounds}{Theorems \ref{thm.bndheight} and \ref{thm.bndwidth}}
\newtheorem*{cor_loceuler}{Corollary \ref{cor.loceuler}}
\theoremstyle{definition}
\newtheorem{defn}[theorem]{Definition}
\newtheorem{exa}[theorem]{Example}
\newtheorem{rem}[theorem]{Remark}

\newcommand{\low}{^{\vphantom{'}}}
\newcommand{\mat}[4]{\bigl(\!\begin{smallmatrix} #1&#2\\#3&#4\end{smallmatrix}\!\bigr)}
\newcommand{\bmat}[4]{\Bigl(\!\begin{smallmatrix} #1&#2\\#3&#4\end{smallmatrix}\!\Bigr)}
\newcommand{\ce}{\mathrel{\mathop:}=}
\newcommand{\ec}{=\mathrel{\mathop:}}
\newcommand{\1}{^{-1}}
\newcommand{\rest}[1]{\vert_{#1}}
\newcommand{\Oh}{\mathcal{O}}
\newcommand{\sM}{\mathcal{M}}
\newcommand{\C}{\mathbb{C}}
\newcommand{\PP}{\mathbb{P}}
\newcommand{\al}{\alpha}
\newcommand{\be}{\beta} 
\newcommand{\ga}{\gamma}
\newcommand{\de}{\delta}
\newcommand{\ze}{\zeta}
\newcommand{\Chat}{{\widehat\ell}}
\newcommand{\iso}{\cong}
\newcommand{\tensor}{\otimes}
\newcommand{\bh}{\mathbf{h}}
\newcommand{\bw}{\mathbf{w}}
\newcommand{\br}{\mathbf{r}}
\newcommand{\bt}{\mathbf{t}}

\newcommand{\op}{\mathcal{O}_{\mathbb{P}^1}}
\DeclareMathOperator{\Ext}{Ext}
\DeclareMathOperator{\Elm}{Elm}

\DeclareMathOperator{\SHom}{\mathscr{H}{\!\mathit{om}}}
\DeclareMathOperator{\SEnd}{\mathscr{E}{\!\mathit{nd}}}
\DeclareMathOperator{\End}{End}
\DeclareMathOperator{\length}{length}
\DeclareMathOperator{\coker}{coker}
\DeclareMathOperator{\Tot}{Tot}

\begin{document}

\maketitle

\addcontentsline{toc}{section}{Abstract}\begin{abstract}\noindent
We study moduli of vector bundles on a two-dimensional neighbourhood
$Z_k$ of an irreducible curve $\ell \iso \PP^1$ with $\ell^2 = -k$ and
give an explicit construction of their moduli stacks.  For the case of
instanton bundles, we stratify the stacks and construct moduli
spaces. We give sharp bounds for the local holomorphic Euler
characteristic of bundles on $Z_k$ and prove existence of families of
bundles with prescribed numerical invariants.  Our numerical
calculations are performed using a \emph{Macaulay~2} algorithm, which
is available for download at \url{http://www.maths.ed.ac.uk/~s0571100/Instanton/}.
\end{abstract}

\section{Introduction}

We study moduli spaces of rank-$2$ bundles on a two-dimensional
neighbourhood of an irreducible curve $\ell \iso \PP^1$ with negative
self-intersection $\ell^2 = -k \neq 0$. We are interested in the
behaviour of bundles over a small analytic neighbourhood of $\ell$
inside a smooth surface $Z_k$, and in coherent sheaves near the
singular point of the surface $X_k$ obtained from $Z_k$ by contracting
the curve $\ell$. For this ``local'' problem of bundles near $\ell$ it
is enough to focus on vector bundles over the total space of
$\Oh_{\PP^1}(-k)$. Hence we take $Z_k \ce
\Tot\bigl(\Oh_{\PP^1}(-k)\bigr)$, where $\ell \subset Z_k$ is the zero
section. We write $\pi \colon Z_k \to X_k$ for the map that contracts
$\ell$ to a point. We give an explicit construction of the moduli 
stack of  rank-$2$ bundles on $Z_k$. 
 
A bundle $E$ over $Z_k$ has \emph{splitting type} $(j_1, \dotsc, j_r)$
if $E\rest{\ell} \iso \bigoplus_{i=1}^r \Oh_{\PP^1}(j_i)$ with $j_1 \leq
\cdots \leq j_r$. For the moduli problem we concentrate on the case of
rank-$2$ bundles $E$ with vanishing first Chern class; in this case the
splitting type of $E$ must be $(-j,j)$ for some $j \geq 0$, and for
short we say that $E$ has splitting type $j$. We then define
 \[ \sM_j(k) = \bigl\{ E\to Z_k : E\rest{\ell} \iso \Oh(j)\oplus\Oh(-j) \bigr\} \bigm/ \sim \]
for the moduli (stack) of bundles over $Z_k$ of splitting type $j$. We
prove:

\begin{thm_genopset}
For $j \geq k$, $\sM_j(k)$ has an open, dense subspace homeomorphic to
a complex projective space \ $\PP^{2j-k-2}$ minus a closed subvariety
of codimension at least $k+1$.
\end{thm_genopset}

The moduli $\sM_j(k)$ contains only one point if $2j - 2 < k$, and it
is non-Hausdorff for $2j - 2 \geq k$. The cases when $j = nk$, that
is, when the splitting type is a multiple of $k = -\ell^2$, are of
special interest for applications to physics, because they correspond
to instantons. We cite:

\begin{theorem*}[{\cite[Corollary 5.5]{GKM}}]
An $\mathfrak{sl}(2,\mathbb{C})$-bundle 
over $Z_k$ represents an instanton if and only if its splitting type 
is a multiple of $k.$
\end{theorem*}

We give a stratification of the instanton moduli stacks, that is of
$\sM_j(k)$ for the case $j = nk$, dividing it into Hausdorff components
(in the analytic topology). For such a stratification, we need
numerical invariants, which we now define: Given a vector bundle $E$
over $Z_k$, we define the Artinian sheaf $Q_E$ on $X_k$ by the exact
sequence
\begin{equation}\label{Q}
  0 \longrightarrow \pi_*E \longrightarrow (\pi_*E)^{\vee\vee}
  \longrightarrow Q_E \longrightarrow 0 \text{ .}
\end{equation}
$E$ has two independent invariants, called \emph{height} and the
\emph{width}. The terminology comes from instantons, cf.\ \cite{NI}.

\begin{defn}\label{def.wh}
We define the \emph{height} and \emph{width} of $E$ by
\[ \bh_k(E) \ce \length R^1 \pi_*E \quad\text{and}\quad \bw_k(E) \ce \length Q_E \text{ .} \]
\end{defn}

We show:

\begin{thm_hdcomps} 
If $j=nk$ for some $n \in \mathbb{N}$, then the pair $(\bh_k,\bw_k)$
stratifies instanton moduli stacks $\sM_j(k)$ into Hausdorff
components.
\end{thm_hdcomps}

\begin{rem}
Note that the invariants are defined for any bundle (or sheaf) over
$Z_k$ but we only claim the stratification result for instantons, that
is, for the case when $j = nk$ is a multiple of $k$. In fact, this is
a necessary condition, and $\sM_3(2)$ is already an example that
justifies the necessity of the ``instanton'' condition (see Example
\ref{ex.noninst}), and application of the embedding theorem
\ref{embed} provides infinitely many such cases.
\end{rem}

One could also decompose $\sM_j(k)$ by the local holomorphic Euler
characteristic $\chi(\ell,E)$ (see Definition \ref{def.euler}).
However, $\sM_j(k)$ has non-Hausdorff subspaces with fixed
$\chi(\ell,E)$, of which the simplest example is $\sM_3(1)$ (see
Example \ref{ex.nonhaus}). See \cite{PAMS} for the case of an
exceptional curve, i.e.\ $k=1$. We prove a simple general formula for
the height:

\begin{thm_height_nonsplit}
Let $E$ be the holomorphic, rank-$2$, non-split vector bundle of
splitting type $j$ which is represented in canonical form \eqref{poly}
by $p$, and let $m > 0$ be the smallest exponent of $u$ appearing in
$p$. With $\mu = \min\bigl(m, \lfloor\frac{j-2}{k}\rfloor\bigr)$, we
have
\[ l(R^1\pi_*E) \geq \mu \left(j - 1  - k\;\frac{\mu - 1}{2}\right) \text{ ,} \]
and equality holds if $p$ is holomorphic on $Z_k$.
\end{thm_height_nonsplit}

We then calculate explicitly sharp bounds for the invariants $(\bh_k, \bw_k)$:

\begin{thm_bounds}
Let $E$ be a holomorphic rank-$2$ vector bundle over $Z_k$ with
$c_1=0$ and splitting type $j > 0$. Let $n_1 = \bigl\lfloor
\frac{j-2}{k} \bigr\rfloor$ and $n_2 = \bigl\lfloor \frac{j}{k}
\bigr\rfloor$. Then the following bounds are sharp:
\begin{gather*}
  j-1 \leq \bh_k(E) \leq (j-1)(n_1+1) - k n_1 (n_1+1)/2 \text{ ,} \\
  0 \leq \bw_k(E) \leq (j+1)n_2 - k n_2 (n_2+1)/2 \text{ ,}
\end{gather*}
and $\bw_1(E) \geq 1$.
\end{thm_bounds}

\begin{defn}[{\cite[Def.\ 3.9]{Bl}}]\label{def.euler}
Let $\sigma(\widetilde{X}, \ell) \to (X,x)$ be a resolution of an
isolated quotient singularity. Let $\widetilde{\cal F}$ be a reflexive
sheaf of rank $n$ on $\widetilde {X}$, and set $\mathcal{F} \ce
(\sigma_* \widetilde{\mathcal{F}})^{\vee\vee}$; notice that there is a
natural injection $\sigma_*\widetilde{\cal F} \hookrightarrow {\cal
F}$. Then the \emph{local holomorphic Euler characteristic} of
$\widetilde{\mathcal{F}}$ is
\[ \chi\bigl(x, \widetilde{\cal F}\bigr) \ce \chi\bigl( \ell, \widetilde{\cal F} \bigr)
    \ce h^0 \bigl(X; \mathcal{F} \bigl/ \sigma_* \widetilde{\cal F}\bigr)
    + \sum_{i=1}^{n} (-1)^{i-1} h^0\bigl(X; R^i \sigma_*\widetilde{\cal F} \bigr) \text{ .} \]
\end{defn}

\begin{cor_loceuler}
Let $E$ be a rank-$2$ bundle over $Z_k$ of splitting type $j > 0$ and
let $j=nk+b$ such that $0 \leq b < k$. The following are sharp bounds
for the local holomorphic Euler characteristic of $E$:
\[ j - 1 \leq \chi( \ell, E) \leq \begin{cases} n^2k+2nb+b-1 &
   \text{if } k \geq 2 \text{ and } 1 \leq b < k \text{ ,} \\
   n^2k & \text{if } k \geq 2 \text{ and } b = 0 \text{ ,} \end{cases} \]
and
\[ j \leq \chi( \ell, E ) \leq j^2 \text{ for } k=1 \text{ .} \]
\end{cor_loceuler}

Next, we consider the question of existence of vector bundles. We
recall the concept of an \emph{admissible sequence} (Definition
\ref{dfn.admissible}) and prove the following existence result:

\begin{thm_admissible}
Fix an admissible sequence $\bigl\{a(i,l)\bigr\}_{i=1}^{t}$ and let
$E$ and $F$ be rank-$r$ vector bundles on $\Chat$ with
$\bigl\{a(i,l)\bigr\}_{i=1}^{t}$ as an associated admissible
sequence. Then there exists a flat family $\bigl\{E_s\bigr\}_{s\in T}$
of rank-$r$ vector bundles on $\Chat$ parametrised by an integral
variety $T$ and $s_0, s_1 \in T$ with $E_{s_0} \iso E$ and $E_{s_1}
\iso F$ such that $E_s$ has $\bigl\{a(i,l)\bigr\}_{i=1}^{t}$ as
admissible sequence for every $s \in T$.
\end{thm_admissible}

In Section~\ref{sec.bounds} we calculate numerical invariants for
bundles near negative curves. In Section~\ref{sec.balancing} we
describe the method of \emph{balancing} bundles and show the existence
of families of bundles with prescribed numerical invariants. In
Section~\ref{sec.moduli} we study moduli of rank-$2$ bundles on $Z_k$.

Our calculations were often performed on a computer using an
implementation of the height and width computations as described in
\cite{GKM} written with the computer algebra software
\emph{Macaulay~2} \cite{M2}. The program can be downloaded from
\url{http://www.maths.ed.ac.uk/} \url{~s0571100/Instanton/}.

The computer program was essential for discovering the results of
Sections \ref{sec.bounds} and \ref{sec.moduli}; algebraic calculations
are straightforward but too laborious to carry out by hand. The
original program was developed by I.\ Swanson in \cite{GS} for bundles
on $Z_1$, where the contraction $\pi \colon Z_1 \to X_1$ has the
lifting map $\widetilde\pi \colon \mathcal{O}_{X_1} \to \pi_*
\mathcal{O}_{Z_1}$ given in coordinates by $x \mapsto u, y \mapsto
zu$. We generalised the computation to the case $Z_k \to X_k$, where
$X_k$ has the coordinate ring $\C[x_0,x_1, \dotsc, x_k] \bigl/\{x_i
x_{i+h} - x_{i+1} x_{i+h-1} \}$ for $0 \leq i \leq i+h \leq k$ and the
lifting map is $x_i \mapsto z^i u$, but the algorithms for the
computation of the module of sections and the width are otherwise
essentially the same.

\paragraph{Acknowledgments.} We would like to thank the editor for 
having done such a wonderful job, and we thank the referee for several
helpful suggestions.

\section{Bounds}\label{sec.bounds}

Let $Z_k$ be the total space of $\Oh_{\PP^1}(-k)$ and $\ell \iso
\PP^1$ the zero section, defined by the ideal sheaf
$\mathcal{I}_\ell$, so that $\ell^2 = -k$. We write
\[ \ell_N = \bigl(\ell, \ \mathcal{O}_{Z_k}\bigl/\mathcal{I}_\ell^{N+1}\rvert_\ell\bigr) \]
for the $N^\mathrm{th}$ infinitesimal neighbourhood of $\ell$, $\Chat
= \varprojlim\ell_N$ for the formal neighbourhood of $\ell$ in $Z_k$,
and $\Oh(j)$ for the line bundle on $Z_k$ or on $\Chat$ that restricts
to $\Oh_{\PP^1}(j)$ on $\ell$. A vector bundle $E$ has \emph{splitting
type} $(j_1, \dotsc, j_r)$ if $E\rest{\ell} \iso \bigoplus_{i=1}^r
\Oh_{\PP^1}(j_i)$ with $j_1 \geq \dotsb \geq j_r$. A result of
Griffiths \cite{GR} implies that $E$ splits on $\Chat$ if $j_1 - j_r
\leq k+1$.

\begin{lemma}\label{nbd}
Let $Z$ be a smooth surface containing a curve $\ell \iso \PP^1$ with
$\ell^2 = -k$. Let $E$ be a rank-$r$ vector bundle on $Z$ of splitting
type $j_1 \geq j_2 \ge \dotsb \geq j_r$, and assume that $j_1 - j_r
\leq k+1$. Then $E$ splits on the formal neighbourhood $\Chat$ of $\ell$,
that is, $E\rest{\Chat} \iso \bigoplus^r_{i-1} \Oh_{\Chat}(j_i)$.
\end{lemma}
\begin{proof}
This follows from \cite[Propositions 1.1 and 1.4]{GR}: The point is
that $\SHom_{\Oh_\Chat}(E, E)\rest{\ell}$ is a direct sum of line bundles of degree
$\geq -(j_1-j_r)$, so the first-order infinitesimal extensions
$H^1\bigl(\Chat; \SHom_{\Oh_\Chat}(E, E) \tensor \Oh(-\ell)\bigr)$ vanish.
\end{proof}

As an aside, we note a relation with the McKay correspondence (see
e.g.\ \cite{Re}), which relates the (co)homology of a resolution
$Y\to\mathbb{C}^2\bigl/G$ to the irreducible representations of a
finite group $G$ of automorphisms: It is known that a reflexive sheaf
on a surface quotient singularity $\C^2 \bigl/ G$ is a direct sum of
the tautological sheaves obtained from the irreducible representations
of $G$. For cyclic quotient singularities $\frac{1}{r}(1,a)$ these are
just the eigensheaves $\Oh(i)$ of the group action, such that
$\pi_*\Oh = \bigoplus_{i=0}^{r-1} \Oh(i)$.

We now study the holomorphic invariants $\bh_k(E) = l(R^1\pi_*E)$ and
$\bw_k(E) = l(Q_E)$ as in Definition \ref{def.wh}.

\begin{rem}\label{TFF}
In principle, we need the Theorem on Formal Functions
\cite[p.~276]{HA}, to calculate $\bw_k$ and $\bh_k$:
\[ \bw_k(E) = \dim_\mathbb{C}\coker\bigl(\rho \colon M \hookrightarrow M^{\vee\vee}\bigr)
   \text{ , \ \ where } M \ce \varprojlim H^0\bigl(\ell_n; E\rvert_{\ell_n}\bigr) \text{ , and} \]
\[ \bh_k(E) = \dim_\mathbb{C} \varprojlim H^1\bigl(\ell_n; E\rvert_{\ell_n}\bigr)
    \text{ .} \]
However, since holomorphic bundles on $Z_k$ are algebraic, the limit
stabilises at a finite order and it is enough to compute the
cohomology on a fixed infinitesimal neighbourhood $\ell_N$, where $N$
is not too small (see \cite[Lemmas 2.1--2.3]{CA2}) --- we restate
these results in a slightly generalised form in Lemma
\ref{lem.finiteN} below.
\end{rem}

It follows from $H^*\bigl(\bigoplus_{i=1}^N \mathcal{F}_i\bigr) =
\bigoplus_{i=1}^N H^*\bigl(\mathcal{F}_i\bigr)$ for coherent sheaves
$\mathcal{F}_i$ that heights and widths behave additively for bundles
that split on the formal neighbourhood $\Chat$:

\begin{proposition}\label{split}
Let $E$ be a holomorphic vector bundle on $Z_k$. If $E\rest{\Chat}
\iso \bigoplus^r_{i=1} \Oh(j_i)$, then
\[
  \bw_k(E) = \sum_{i=1}^r \bw_k\bigl(\Oh(j_i)\bigr) \quad\hbox{and}\quad
  \bh_k(E) = \sum_{i=1}^r \bh_k\bigl(\Oh(j_i)\bigr) \text{ .} \qquad\qed 
\]
\end{proposition}

\begin{corollary}
Let $E$ be a holomorphic rank-$r$ vector bundle over $Z_k$ such that
$E\rest\ell \iso \bigoplus^r_{i=1} \Oh_\ell(j_i)$, with $j_1 \geq j_2 \geq
\dotsb \geq j_r$ and $j_r -j_1 \geq -k-1$. Then
\[
  \bw_k(E) = \sum_{i=1}^r \bw_k\bigl(\Oh(j_i)\bigr) \quad\hbox{and}\quad
  \bh_k(E) = \sum_{i=1}^r \bh_k\bigl(\Oh(j_i)\bigr) \text{ .}
\]
\end{corollary}
\begin{proof}
By Lemma~\ref{nbd}, $E\rest{\Chat} \iso \bigoplus^r_{i=1}
\Oh_{\Chat}(j_i)$.  The result follows from Proposition~\ref{split}
(and hence $E$ has the same invariants as the split bundle).
\end{proof}

We fix once and for all coordinate charts on $Z_k$, to which we will
refer as \emph{canonical coordinates},
\begin{equation}\label{eq.cancoords}
  U = \C^2_{z,u} = \bigl\{z,u\bigr\} \qquad\text{and}\qquad
  V = \C^2_{\ze,v} = \bigl\{\ze, v\bigr\} \text{ ,}
\end{equation}
glued by $\ze =z\1$ and $v=z^ku$. In these charts the bundle $\Oh(j)$
has the transition matrix $\bigl(z^{-j}\bigr)$.

Bundles with vanishing first Chern class are of special interest for
applications to physics. For example, the bundles with splitting type
$nk$ correspond to framed instantons under a local version of the
Kobayashi--Hitchin correspondence (see \cite{GKM} and \cite{LT}).  For
the case $c_1=0$, we also calculate the lower bounds for the numerical
invariants $\bh$ and $\bw$. The second author proved in \cite{CA1} that
holomorphic bundles on $Z_k$ are algebraic extensions of line
bundles. By \cite[Theorem 3.3]{CA1}, a bundle $E$ that is an extension
\begin{equation}\label{eq!ext}
  0 \longrightarrow \Oh(j_1) \longrightarrow E \longrightarrow \Oh(j_2) \longrightarrow 0
\end{equation}
(with $j_1 \leq j_2$) has transition matrix $T =
\mat{z^{-j_1}}{p(z,u)}{0}{z^{-j_2}}$ in canonical coordinates
\eqref{eq.cancoords}, where the extension class
\[ \bigl[p\bigr] \in \Ext^1_{\Oh_{Z_k}}\!\bigl(\Oh(j_2),\; \Oh(j_1)\bigr) \]
may be represented by a \emph{canonical form}
\begin{equation}\label{poly}
  p(z,u) = \sum_{r=1}^{\left\lfloor \vphantom{\textstyle X}(j_2 - j_1 - 2)/k \right\rfloor}
  \sum_{s = kr + j_1 + 1}^{j_2 - 1} p_{rs} \, z^s u^r \text{ .}
\end{equation}

\begin{lemma}[Computation on finite neighbourhoods]\label{lem.finiteN}
Let $E$ be a rank-$2$ bundle of splitting type $j$ on $Z_k$ and set $N \ce
\left\lfloor (2j-2)/k\right\rfloor$. Then the height and width of $E$ are
determined already on $\ell_N$ in the following sense:
\begin{itemize}
\item $\bh_k(E) = \dim_\mathbb{C} H^1\bigl(\ell_N; E\rvert_{\ell_N}\bigr)$.
\item Set $M = H^0\bigl(\ell_N; E\rvert_{\ell_N}\bigr)$ and let
      $\rho\colon M \hookrightarrow M^{\vee\vee}$ be the evaluation
      map. Then $\bw_k(E) = \dim_\mathbb{C} \coker\rho$.
\end{itemize}
\end{lemma}
\begin{proof}
The crucial fact is that the extension class is given by a polynomial
$p$ whose degree in $u$ is at most $N \ce \left\lfloor
(2j-2)/k\right\rfloor$. Now it was shown in \cite[Lemmas 2.1]{CA2} for
the case $k=1$, but readily generalised to all $k \geq 1$, that for $n
\geq N$, the modules $H^i\bigl(\ell_n; E\rvert_{\ell_n}\bigr)$ and
$H^i\bigl(\ell_N; E\rvert_{\ell_N}\bigr)$ have the same module
structure as $\mathcal{O}_0^{\wedge}$-modules for all $i\geq0$, where
$\pi(\ell) = 0 \in X_k$. Hence by \cite[p.~193]{HA}, the inverse
limits from Remark \ref{TFF} stabilise at $n=N$, which proves the
Lemma.
\end{proof}

\subsection{Heights}

We begin with the computation of the height of a line bundle:

\begin{lemma}\label{lem.height_split}
Assume $j \geq 0$ and let $n_1 = \left\lfloor \frac{j-2}{k}
\right\rfloor$. Then
\[ \bh_k\bigl(\Oh(-j)\bigr)= \begin{cases} \displaystyle (j-1)(n_1+1) - kn_1(n_1+1)/2 &
   \text{if $j\geq 2$,} \\ 0 & \text{otherwise.} \end{cases} 
\]
\end{lemma}
\begin{proof}
$R^1\pi_* \Oh(j) = \varprojlim H^1\bigl(\Oh_{\ell_n}(j)\bigr)$, with
surjective restriction maps. The result comes from the exact sequences
\[ 0 \longrightarrow H^1\bigl(\Oh_\ell(-n\ell) \tensor \Oh(j)\bigr)
   \longrightarrow H^1\bigl(\Oh_{\ell_{n+1}}(j)\bigr) \longrightarrow
   H^1\bigl(\Oh_{\ell_n}(j)\bigr) \longrightarrow 0 \]
together with
\[ H^1\bigl(\Oh_\ell(-n\ell) \tensor \Oh(j)\bigr) =
   H^1\bigl(\PP^1; \Oh(j+nk)\bigr) \text{ ,} \]
which gives $\bh_k\bigl(\Oh(-j)\bigr) = \sum_{n=0}^\infty (j-1-nk)^+$,
where $^+$ means the sum of positive terms only.
\end{proof}

This result together with Proposition \ref{split} allows us to compute
the heights of split bundles. For non-split bundles with $c_1=0$ and
splitting type $j \ce j_2 = -j_1$, we use the following:

\begin{theorem}\label{thm.height_nonsplit}
Let $E$ be the non-split bundle of splitting type $j$ represented in
canonical form \eqref{poly} by $p$, and let $m > 0$ be the smallest
exponent of $u$ appearing\footnote{A rank-$2$ bundle whose extension
class is a polynomial not divisible by $u$ is in fact not of splitting
type $j$, but of a lower splitting type. See Remark
\ref{rem.udivides}.} in $p$. With $\mu = \min\bigl(m,
\lfloor\frac{j-2}{k}\rfloor\bigr)$, we have
\[ l(R^1\pi_*E) \geq \mu \left(j - 1  - k\;\frac{\mu - 1}{2}\right) \text{ ,}  \]
and equality holds if $p$ is holomorphic on $Z_k$.
\end{theorem}
\begin{proof}
By Lemma~\ref{lem.h1rep} below, a cocycle in $H^1(E)$ has the
canonical representation on the $U$-chart
\[ \sum_{r=0}^{\left\lfloor\frac{j-2}{k}\right\rfloor}
   \sum_{s=kr-j+1}^{-1} \begin{pmatrix}a_{rs}\\0\end{pmatrix} z^s u^r \text{ .} \]
In this representation, every monomial term $(a_{rs}z^s u^r, 0)$ with
$r<m$ represents a non-trivial cocycle by Lemma
\ref{lem.h1nontriv}. Lastly, if $p$ is holomorphic in $Z_k$, it has
only terms $p_{rs} z^s u^r$ with $0\leq s \leq kr$, and then all terms
with $r \geq m$ are coboundaries by Lemma \ref{lem.h1triv}.
\end{proof}

Putting the split and the non-split cases together, we get sharp
bounds on the heights:

\begin{theorem}\label{thm.bndheight}
Let $E$ be a rank-$2$ bundle over $Z_k$ of splitting type $j > 0$.
Set $n_1 = \left\lfloor \frac{j-2}{k} \right\rfloor$. The following
bounds are sharp:
\[ j-1 \leq \bh_k(E) \leq (j-1)(n_1+1)-k (n_1+1)n_1/2 \text{ .} \]
\end{theorem}
\begin{proof}
The upper bound is attained by the split bundle by Proposition
\ref{prop.splitmax}. In this case, apply Lemma~\ref{lem.height_split}
and Proposition~\ref{split}.

For the lower bound, note that the expression in
Theorem~\ref{thm.height_nonsplit} is always less than in or equal to
the split case. The global minimum is $j-1$, attained with $\mu=1$ in
Theorem~\ref{thm.height_nonsplit}, attained by $p(z,u)=zu$.
\end{proof}

We finish this subsection by proving the details of Theorem
\ref{thm.height_nonsplit}. Recall once and for all that in the charts
$Z_k = U \cup V$ given by \eqref{eq.cancoords}, a function is
holomorphic on $U$ if it is holomorphic in $\{z,u\}$, and on $V$ if it
is holomorphic in $\{z^{-1}, z^ku\}$.

\begin{lemma}\label{lem.h1rep}
Every $1$-cocycle in $H^1(E)$ has a representative of the form
\[ \sum_{r=0}^{\left\lfloor \frac{j-2}{k} \right\rfloor} \sum_{s=kr-j+1}^{-1}
   \begin{pmatrix}a_{rs} \\ 0\end{pmatrix} z^s u^r \text{ ,} \]
with $a_{rs} \in \C$.  In particular, every $1$-cochain
represented by $\left(\begin{smallmatrix}a_{rs} \\ 0\end{smallmatrix}\right) z^s u^r$ with $r, s \geq 0$ is a coboundary.
\end{lemma}
\begin{proof}
Let $\sigma$ be a $1$-cocycle and let $\sim$ denote cohomological
equivalence. A power series representative for a $1$-cochain has the
form
\[ \sigma = \sum_{r=0}^\infty \sum_{s=-\infty}^{\infty} \begin{pmatrix}a_{rs} \\ b_{rs} \end{pmatrix} z^s u^r \text{ ,} \]
with $a_{rs}, b_{rs} \in \C$. The $1$-cochain $s_1 =
\sum_{r=0}^\infty \sum_{s=0}^\infty \bigl( \begin{smallmatrix}a_{rs}
\\ b_{rs}\end{smallmatrix} \bigr) z^s u^r$ is holomorphic on $U$,
hence represents  a coboundary. Consequently
\[ \sigma \sim \sigma - s_1 = \sum_{r=0}^{\infty} \sum_{s=-\infty}^{-1} \begin{pmatrix} a_{rs} \\ b_{rs} \end{pmatrix} z^s u^r \text{ .} \]
Now let $T=\begin{pmatrix} z^j & p \\ 0 & z^{-j} \end{pmatrix}$ be the
transition function of $E$, so that after a change of coordinates,
\[ T \sigma = \sum_{r=0}^{\infty} \sum_{s=-\infty}^{-1} \begin{pmatrix}
   z^j a_{rs} + p\,b_{rs} \\ z^{-j} b_{rs}\end{pmatrix} z^s u^r \text{ .} \]
However, given that $s_2 = \sum_{r=0}^{\infty} \sum_{s=-\infty}^{-1}
\bigl( \begin{smallmatrix} 0 \\ z^{-j} b_{rs} \end{smallmatrix} \bigr)
z^s u^r$ is holomorphic on $V$,
\[ T \sigma \sim T \sigma - s_2 = \sum_{r=0}^{\infty} \sum_{s=-\infty}^{-1}
   \begin{pmatrix} z^j a_{rs} + p\,b_{rs} \\ 0 \end{pmatrix} z^s u^r \text{ ,} \]
and going back to the $U$-coordinate chart,
\[ \sigma = T^{-1} T \sigma \sim \sum_{r=0}^{\infty} \sum_{s=-\infty}^{-1}
   \begin{pmatrix} a_{rs} + z^{-j} p\,b_{rs} \\ 0 \end{pmatrix} z^s u^r \text{ .} \]
But $p$ contains only terms $z^k$ for $k \leq j-1$, therefore $z^{-j}
p$ contains only negative powers of $z$. Renaming the coefficients we
may write
\[ \sigma = \sum_{r=0}^{\infty} \sum_{s=-\infty}^{-1} \begin{pmatrix} a'_{rs} \\ 0 \end{pmatrix} z^s u^r \]
for some $a'_{rs} \in \C$, and consequently $T\sigma =
\sum_{r=0}^{\infty} \sum_{s=-\infty}^{-1} \bigl( \begin{smallmatrix}
z^j a'_{rs} \\ 0 \end{smallmatrix} \bigr) z^s u^r$. Here each term
$a'_{rs}z^{j+s} u^r$ satisfying $j+s \leq kr$ is holomorphic on the
$V$-chart. Subtracting these holomorphic terms we are left with an
expression for $a$ where the index $s$ varies as $kr-j+1 \leq s \leq
-1$. This in turn forces $r \leq \bigl\lfloor \frac{j-2}{k}
\bigr\rfloor$, giving the claimed expression for the $1$-cocycle.
\end{proof}

\begin{lemma}\label{lem.h1nontriv}
Let $E$ be the non-split bundle over $Z_k$ represented in canonical
form by $(j,p)$ and $m > 0$ the smallest exponent of $u$ appearing in
$p$. If $\mu = \min\bigl(m, \lfloor\frac{j-2}{k}\rfloor\bigr)$, then
\[ l(R^1\pi_*E) \geq \mu\left(j - 1 - k\;\frac{\mu-1}{2}\right) \text{ .} \]
\end{lemma}
\begin{proof}
Assume $m \leq \bigl\lfloor \frac{j-2}{k} \bigr\rfloor$. Let $\sigma =
(a,0)$ denote a $1$-cocycle where $a = z^s u^r$ with $0 \leq r \leq
m-1$ and $kr - j + 1 \leq s \leq - 1$ (due to Lemma~\ref{lem.h1rep}).
We claim that $\sigma$ represents a non-zero cohomology class. In
fact, since $\sigma$ is not holomorphic on $U$, for $\sigma$ to be a
coboundary there must exist a $U$-coboundary $\alpha$ such that
$\sigma+\alpha$ is holomorphic on $V$. Hence, there must exist a
polynomial $X$, holomorphic on $U$, making the expression $z^j a + p
X$ holomorphic on $V$. However, $z^j a = z^{s+j}u^r$ is not
holomorphic on $V$, since $s+j \geq kr+1$. Moreover, by the choice of
$m$, no term in $p X$ cancels $z^j a$. Consequently, no choice of $X$
solves the problem of holomorphicity on $V$. Hence $l(R^1\pi_*E)$ is
at least the number of independent cocycles of the form $\sigma =
(a,0)$, where $a = z^s u^r$ with $0 \leq r \leq m-1$ and $kr - j + 1
\leq s \leq -1$. But by exactly the same reasoning, any linear
combination $p_{rs} z^s u^r$ (with $r$ and $s$ satisfying the above
inequalities) is a coboundary if and only if every term is a
coboundary individually, and hence the monomial cocycles $\sigma =
(a,0)$ are linearly independent. There are $m\bigl((j-1) -
k(m-1)/2\bigr)$ such terms.

On the other hand, if $m > \bigl\lfloor \frac{j-2}{k} \bigr\rfloor$,
then by an analogous argument, none of the monomial terms in the
canonical cocycle form of Lemma~\ref{lem.h1rep} are coboundaries, and
so the same formula holds with $m$ replaced by $\bigl\lfloor
\frac{j-2}{k} \bigr\rfloor$.
\end{proof}

\begin{lemma}\label{lem.h1triv}
Let $E$ satisfy the same conditions as in Lemma~\ref{lem.h1nontriv},
and assume in addition that $p$ is holomorphic on $Z_k$.
Then
\[ l(R^1\pi_*E) = \mu\left(j - 1 - k\;\frac{\mu-1}{2}\right) \text{ .} \]
\end{lemma}
\begin{proof}
If $p$ is holomorphic on $Z_k$, it can be written as
\[ p(z,u) = \sum_{r=1}^{\left\lfloor (2j-2)/k \right\rfloor} \ \sum_{s = 0}^{\min(j-1,\;kr)} p_{rs} z^s u^r \text{ .} \]

Following the proof of Lemma \ref{lem.h1nontriv}, we now need to show
that all remaining $1$-cocycles, namely $\sigma =
\bigl(a_{\bar{r}\bar{s}} z^{\bar{s}} u^{\bar{r}}, 0\bigr)$ with
\[ m \leq \bar{r} \leq \left\lfloor\frac{j-2}{k}\right\rfloor \text{ \ and \ }
   k\bar{r} - j + 1 \leq \bar{s} \leq -1 \tag{$*$} \]
are coboundaries. A fixed such $1$-cocycle $\sigma = \bigl(a z^{\bar{s}}
u^{\bar{r}}, 0\big)$ is a $1$-coboundary if and only if there exist
$0$-cocycles $(\alpha, \beta) \in H^0(V; E)$ and $(x,y) \in H^0(U; E)$
such that
\[ T^{-1} \begin{pmatrix}\alpha \\ \beta \end{pmatrix} + \sigma + \begin{pmatrix}x \\ y \end{pmatrix} = 0 \text{ .} \]
In other words, we need to find $\alpha$, $\beta$ holomorphic on $V$
such that
\[ A = \begin{pmatrix}z^{-j} & -p \\ 0 & z^{j} \end{pmatrix}
   \begin{pmatrix} \alpha \\\beta \end{pmatrix} + \begin{pmatrix} a z^{\bar{s}} u^{\bar{r}} \\ 0 \end{pmatrix} = 
   \begin{pmatrix} z^{-j} \alpha -p \beta + a z^{\bar s} u^{\bar r} \\  z^{j}\beta \end{pmatrix} \]
is holomorphic on $U$. Our goal is to cancel the coefficient $a$.
Notice that, since $z^j az^{\bar s} u^{\bar r}$ is not holomorphic on
$V$, no choice of $\alpha$ can cancel $a$.  So, we can assume $\alpha
= 0$ and study the problem of holomorphicity on $U$ for the simplified
matrix
\[ A = \begin{pmatrix} -p \beta + az ^{\bar s} u^{\bar r} \\  z^{j}\beta \end{pmatrix} \text{ .} \]
Let $\sup \{ l : p_{ml} \neq 0 \} \ec l_0 \geq 0$ by assumption, and 
set
\[ \beta = \frac{a z^{\bar s} u^{\bar r}}{p_{ml_0} z^{l_0} u^m} = 
   \frac{a}{p_{ml_0}} z^{\bar s - l_0} u^{\bar r-m} \text{ .} \]
Since $\bar s - l_0 <0$, $\beta$ is holomorphic on $V$. Now the matrix
$A$ becomes
\[ A = \begin{pmatrix} - \left(\sum_{l=0}^{l_0-1}
   p_{ml} z^l u^m + p' \right) \frac{a}{p_{ml_0}} z^{\bar s - l_0} u^{\bar r-m} \\
   z^{j} \frac{a}{p_{m l_0}} z^{\bar s - l_0} u^{\bar r-m} \end{pmatrix} \text{ ,} \]
where $u^{m+1}$ divides $p'$. We have $l_0 \leq \min(j-1,\;km)$ and
$m\leq\bar{r}$ by assumption. Using $(*)$, we obtain $j + \bar s - l_0
\geq 0$; consequently the second coordinate of $A$ is holomorphic on
$U$. The first coordinate of $A$ is
\[ -\sum_{l=0}^{l_0-1} \frac{a\,p_{ml}}{p_{ml_0}} z^{l + \bar s - l_0} u^{\bar r}
   - p' \frac{a}{p_{ml_0}} z^{\bar s - l_0} u^{\bar r-m} =
   -\sum_{l=-l_0}^{-1} p_{ml} z^{l+ \bar s} u^{\bar r} - p' \frac{a}{p_{ml_0}} z^{\bar s - l_0} u^{\bar r-m} \text{ .} \]
We obtain $A \sim -A^- - A^+$, where
\[ A^- \ce \sum_{l = -l_0}^{-1} \begin{pmatrix}p_{ml} z^{l + \bar s} u^{\bar r} \\ 0 \end{pmatrix} \text{ ,\quad}
   A^+ \ce \begin{pmatrix}p' \frac{a}{p_{ml_0}} z^{\bar s - l_0} u^{\bar r-m} \\ 0 \end{pmatrix} \text{ .} \]
This expresses $\sigma$ as a linear combination of $1$-cocycles having
distinct degrees. Notice that $A^-$ is either zero, or contains only
powers of $z$ strictly lower than $\bar s$; whereas $A^+$ is either
zero or contains only powers of $u$ strictly greater than $\bar r$.

To complete the proof we now proceed by finite induction on the
exponents $\bar r$ and $\bar s$. The initial case is when $\bar r=m$
and $\bar s=km-j+1$, and in this case it follows from Lemma
\ref{lem.h1rep} that $A^- \sim 0$, and we are left with $A^+$, whose
degree in $u$ is greater than $m$.

Now assume the lemma has been proved for cycles of the form $(a_{rs}
z^s u^r, 0)$ for all $r < \bar r$ and $s < \bar s$ and apply the same
algorithm to $\sigma = (a_{\bar r \bar s} z^{\bar s} u^{\bar r},
0)$. Then again, $A^-$ gets expressed as a combination of cocycles
already known to be null-cohomologous, whereas $A^+$ gets expressed in
terms of cycles of higher degree in $u$.

The induction finishes at $r=\lfloor (j-2)/k\rfloor$ when 
$A^+$ also becomes zero by Lemma \ref{lem.h1rep}.
\end{proof}

Even though we did not compute explicitly the height for bundles whose
extension class $p$ is \emph{not} holomorphic on $Z_k$ (the computer
algorithm can compute it for arbitrary $p$), we can show that the
split bundle attains the maximally possible height for a fixed
splitting type $j$:

\begin{proposition}\label{prop.splitmax}
Let $E(j,p)$ be the bundle of splitting type $j$ whose extension class
is given by $p$, and let $\bar{E}(j) \ce \mathcal{O}(-j) \oplus
\mathcal{O}(j)$ denote the split bundle. If $u|p(z,u)$ and
$p\not\equiv0$, then
\[ \mathbf{h}_k\bigl(\bar{E}(j)\bigr) \geq \mathbf{h}_k\bigl(E(j,p)\bigr) \text{ .} \]
\end{proposition}
\begin{proof}
Applying $\pi_*$ to the short exact sequence 
\begin{equation}\label{eq.SES-E}
  0 \longrightarrow \mathcal{O}(-j) \xrightarrow{ \ \imath\ } E \longrightarrow \mathcal{O}(j) \longrightarrow 0
\end{equation}
gives the long exact sequence
\[ \dotsb \longrightarrow \pi_*{\mathcal O}(j) \xrightarrow{ \ \delta\ } R^1\pi_*\mathcal{O}(-j) 
   \xrightarrow{ \ \imath_*\ } R^1\pi_*E \longrightarrow R^1\pi_*\mathcal{O}(j) \longrightarrow 0 \text{ ,} \]
with $\pi_*{\mathcal O}(j) \neq 0$ because $ j \geq 0$; and the
connecting homomorphism $\delta$ is zero precisely when the sequence
\eqref{eq.SES-E} splits. Hence, if $E$ does not split,
\begin{align*}
  l\bigl(R^1\pi_*E\bigr) &\leq l\bigl(\imath_* R^1 \pi_* \mathcal{O}(-j)\bigr) + l\bigl(R^1\pi_*\mathcal{O}(j)\bigr) \\
  &\leq l\bigl(R^1\pi_*\mathcal{O}(-j)\bigr) + l\bigl(R^1\pi_*\mathcal{O}(j)\bigr) \text{ ,}
\end{align*}
with equality holding for the split bundle (and possibly other bundles
as well, for instance in cases when $k$ does not divide $j$).
\end{proof}

\begin{rem}\label{rem.udivides}
Let us stress again that it is necessary for $u$ to divide $p$ in
order for the bundle $E(j,p)$ to split as $\op(-j)\oplus\op(j)$ on
$\ell = \bigl\{u=0\bigr\}$: For example, the bundle on $Z_1$ given by
transition function $\mat{z^2}{z}{0}{z^{-2}}$, which we would denote
by $E(2,z)$, is in fact not of splitting type $2$, but $1$:
\[ \begin{pmatrix} 1 & 0 \\ -z^{-3} & 1 \end{pmatrix} \begin{pmatrix} z^2 & z \\ 0 & z^{-2} \end{pmatrix}
   \begin{pmatrix} 0 & -1 \\ 1 & z \end{pmatrix} = \begin{pmatrix} z & 0 \\ 0 & z^{-1} \end{pmatrix} \]
\end{rem}

\subsection{Widths}\label{sec.widths}

Again we begin by computing the width of \emph{line} bundles:
Since computations get quite complicated very quickly, we 
chose to present one explicit example, and refer to the 
program for the general case.

\begin{exa}[Computation of $l(Q)$ for the bundle $\mathcal{O}(-3)$ over 
$Z_2$]\label{lQexample}\mbox{}\\
Here we show that ${\bf w}_2\bigl(\Oh(-3)\bigr)=0$:

In our canonical coordinates $Z_2$ is given by charts $U =
\bigl\{(z,u)\bigr\},$ $V = \bigl\{(z^{-1}, z^2u)\bigr\}$, and the
transition matrix for the bundle $E = \Oh(-3)$ is $T =
\bigl(z^3\bigr)$. $X_2$ is the singular space obtained by the
contraction of the zero section $\pi \colon Z_2 \to X_2$. The
coordinate ring of $X_2$ is $\mathbb{C}[x,y,w] \bigl/ (y^2-xw)$, and
the lifting map induced by $\pi$ is (on the $U$-chart) $x \mapsto u$,
$y \mapsto zu$ and $w \mapsto z^2u$.

Let $M = (\pi_* E)_0^\wedge$ denote the completion of the stalk
$(\pi_*E)_0$. Let $\rho \colon M \hookrightarrow M^{\vee\vee}$ denote
the natural inclusion of $M$ into its double dual. We want to compute
$l(Q) = \dim(\coker\rho)$. By the Theorem on Formal Functions,
\[ M \cong \varprojlim H^0\bigl(\ell_n, E\lvert_{\ell_n}\bigr) \text{ ,} \]
where $\ell_n$ denotes the $n^\text{th}$ infinitesimal neighbourhood
of $\ell$. To determine $M$, it suffices by Lemma \ref{lem.finiteN} to
calculate $H^0\bigl(\ell_{n}, E|_{\ell_{n}}\bigr)$ for a fixed
$n \geq N = \left\lfloor(2j-2)/k\right\rfloor$ (cf.\ also \cite{FM} or
\cite{CA1}) and the relations among its generators under the action of
$\Oh_0^{\wedge} \cong \mathbb{C}[[x,y,w]]\bigl/(y^2-xw)$. In this
example, a section $\sigma$ of $E$ on the $U$-chart has the expression
\[ \sigma = \sum_{r=0}^{N} \sum_{s=0}^{\infty} a_{rs} z^r u^s \text{ ,} \]
and changing coordinates, we have the condition that
\[ T\sigma = z^3 \sigma \]
must be holomorphic on $V = \bigl\{(z^{-1}, z^2u)\bigr\}$, and hence
must have only terms of the form $z^s u^r$ that satisfy $2 s \leq r$,
with the remaining coefficients vanishing. It follows that the
expression for $\sigma$ has the form
\begin{align*}
  \sigma &= \sum_{r=2}^N \sum_{s=0}^{\bigl\lfloor \frac{2r-3}{2}\bigr\rfloor} a_{rs}{z^s u^r} \\
  &= a_{20} u^2 + a_{21} z u^2 + a_{30} u^3 + a_{31} z u^3 + a_{32} z^2 u^3 + a_{33} z^3 u^3 + \dotsb \text{ .}
\end{align*}
Consequently, all terms of $\sigma$ are generated over
$\Oh_0^{\wedge}$ by $\beta_0 = u^2 = x^2$ and $\beta_1 = zu^2 =
xy$. We obtain the $\Oh_0^{\wedge}$-module $M \cong
\Oh_0^{\wedge}\bigl[\beta_0, \beta_1\bigr] \bigl /{\bf R^1}$, where
${\bf R^1}$ is the set of relations
\[ \mathbf{R^1} = \begin{cases} R^1_1: & \beta_0 y - \beta_1 x \text{ ,} \\ R^1_2: & \beta_0 w - \beta_1 y \text{ .} \end{cases} \]
Consequently, the dual is $M^\vee = \Oh_0^{\wedge}\bigl[\beta_0^{\vee}, \beta_1^{\vee}\bigr] \bigl/ \mathbf{S}$, where
\[ \beta_0^\vee = \begin{cases} \beta_0 \mapsto x \\ \beta_1 \mapsto y \end{cases} \text{ ,}\qquad
   \beta_1^\vee = \begin{cases} \beta_0 \mapsto y \\ \beta_1 \mapsto w \end{cases} \text{ ,} \]
and $\mathbf{S}$ is the set of relations
\[ \mathbf{S} = \begin{cases} S_1: & \beta_0^\vee y - \beta_1^\vee x \text{ ,} \\ S_2: & \beta_0^\vee w - \beta_1^\vee y \text{ .} \end{cases} \]
Clearly, $M \cong M^{\vee}$, consequently $\rho \colon M
\hookrightarrow M^{\vee\vee}$ is also an isomorphism, so ${\bf
w}_2\bigl(\Oh(-3)\bigr)=0$.
\newline\mbox{} 
\hfill$\bigr/\!\!\bigl/$
\end{exa}

\medskip\noindent We remark as an aside that there is another set of
relations ${\bf R^2}$ among the relations of $M$:
\[ \mathbf{R^2} = \begin{cases} R^2_1: & R^1_1 y - R^1_2 x \text{ ,} \\ R^2_2: & R^1_1 w - R^1_2 y \text{ ;} \end{cases} \]
and so on:
\[ \mathbf{R^n} = \begin{cases} R^n_1: & R^{n-1}_1 y - R^{n-1}_2 x \text{ ,} \\ R^n_2: & R^{n-1}_1 w - R^{n-1}_2 y \text{ .} \end{cases} \]
This is an example of a theorem of Eisenbud \cite{eisenbud1980}, which
says that every minimal resolution of a finitely generated module over
$A \ce R \bigl/ (x)$ becomes periodic with period $1$ or $2$ after at
most $\dim A$ steps, where $R$ is a regular ring and $x$ a non-unit.

\begin{lemma}\label{lem.width_split}
Assume $j \geq 0$ and let $n_2 = \left\lfloor\frac{j}{k}\right\rfloor$. Then
\begin{align*}
  \bw_k\bigl(\Oh(j)\bigr)  &= (j+1)n_2 - kn_2(n_2+1)/2 \text{ ,} \\
  \bw_k\bigl(\Oh(-j)\bigr) &= 0 \text{ .}
\end{align*}
\end{lemma}
\begin{proof}
The length $l(Q)$ equals the dimension of $Q_0^{\wedge}$ as a
$\C$-vector space, where $0 \in X_k$ is the singular point. Since $Q$
is defined by the sequence \eqref{Q}, $Q_0^{\wedge}$ is the cokernel
of the map $\bigl(\pi_*\Oh(j)\bigr)_0^{\wedge} \to
\bigl(\pi_*\Oh(j)\bigr)_0^{\vee\vee\wedge}$. By the Theorem on Formal
Functions \cite[p.~276]{HA},
\[ \bigl(\pi_* \Oh(j)\bigr)_0^{\wedge} = \varprojlim H^0
   \Bigl(\ell_n; \Oh(j)\bigl\lvert_{\ell_n}\Bigr) \ec M \text{ .} \]
But the limit stabilises at a finite stage by Lemma \ref{lem.finiteN},
so it suffices to compute the group $H^0\bigl(\ell_n;
\Oh(j)\rest{\ell_n}\bigr)$ for large $n$. We first compute the
$k_0$-module structure on $M \ce (\pi_* \Oh(j))_0^{\wedge}$. Note that
here
\[ k_0 \iso \C[[x_0,x_1, \dotsc, x_k]] \bigl/ \bigl\{x_i x_j - x_{i+1} x_{j-1} \bigr\}
   \text{ \ for } i = 0, 1, \dotsc, k-2 \text{ , } j = i+2, \dotsc, k \text{ ,} \]
and the contraction map $\pi \colon Z_k \to X_k$ is given in
$(z,u)$-coordinates by $x_i = z^i u$.

In the case of $\Oh(j)$, $M$ is generated as a $k_0$-module by the
monomials $\be_i = z^i$ for $0 \leq i \leq j$ with relations $\be_i
x_{s-1} - \be_{l-1} x_l = 0$ for $1 \leq i \leq j$ and $1 \leq l\leq
k$. Now, one can use any computer algebra package (e.g.\
\emph{Macaulay~2} \cite{M2}, \emph{Maple}, \emph{Magma}, see
\cite{GS}) to compute $M^{\vee\vee} \bigl/ M$: Writing $j=mk+b$ with
$1\leq b \leq k$ (so that $m=\lfloor (j-1)/k \rfloor$), one finds that
$M^\vee$ and $M^{\vee\vee}$ are generated by $k-b+1$ elements and that
the cokernel has, as a $\mathbb{C}$-vector space, a basis consisting
of all the monomials of degree $\leq m$. Accounting for the relations
coming from $M^{\vee\vee}$ and from the ground ring $k_0$ shows that
the dimension of the cokernel is $(j+1)n_2 - kn_2(n_2+1)/2$.

In the case $\Oh(-j)$, set $\nu = -j \mod k$, so that $-j = -q k +
\nu$ with $0 \leq \nu < k$. Then for large $n$, $H^0\bigl(\ell_n;
\Oh(j)\rest{\ell_n}\bigr)$, and hence $M$, is generated by the set of
monomials $\al_i = z^i u^q$, for $0 \leq i \leq \nu$, with relations
$\al_i x_{l-1} - \al_{i-1} x_l = 0$ for $1 \leq i \leq \nu$ and $1
\leq l \leq k$. By direct calculation, or using a computer algebra
package, one shows that the evaluation map $M \to M^{\vee\vee}$ is an
isomorphism, so $\bw_k(\Oh(-j))=0.$ Explicit computation for $\Oh(-2)$
is given in Example \ref{lQexample}.
\end{proof}

\begin{theorem}\label{thm.bndwidth}
Let $E$ be a rank-$2$ bundle over $Z_k$ of splitting type $j$. Then
the following bounds are sharp: For $j>0$ and with $n_2 = \bigl\lfloor
\frac{j}{k} \bigr\rfloor$,
\[ 0 \leq \bw_k(E) \leq (j+1)n_2 - kn_2(n_2+1)/2 \text{ , and } \bw_1(E) \geq 1 \text{ .} \]
Furthermore, for all $0 < j < k$, $\bw_k(E) = 0$ for all bundles $E$ (and necessarily $k>1$).
\end{theorem}
\begin{proof}
The upper bound is realised by the split bundle, so the right-hand
side follows from Lemma~\ref{lem.width_split} and
Proposition~\ref{split}. This follows from direct computation: First
we write down a set of possible generators of $M$, then we compute
their relations. But relations can only come from the off-diagonal
part of the transition function of $E$, i.e.\ from $p$, so the split
bundle with $p=0$ has the fewest possible relations in $M$, hence the
biggest cokernel $M^{\vee\vee}\bigl/M$.

To calculate the lower bound for $\bw_k(E)$, note that by definition
$\bw_k(E) \geq 0$. But for $k>1$, the bundle given in canonical
coordinates by transition matrix $\mat{z^{j}}{zu}{0}{z^{-j}}$ has
width zero. Direct computations as described in \cite{CA2} show that
for the case $k=1$, $\bw_1(E) = 1$. 
\end{proof}

\begin{rem}
It is interesting to notice that on the space $Z_1$ any nontrivial
rank-$2$ bundle $E$ with $c_1(E)=0$ has ${\bf w}_1(E) \neq 0$. This is
in strong contrast with what happens on $Z_k$ for $k>1$.
\end{rem}

\subsection{Local holomorphic Euler characteristic}

By summing up the results for heights and widths and using
Definition~\ref{def.euler} we obtain bounds on the local Euler
characteristic of a rank-$2$ bundle $E$ on $Z_k$ near $\ell$. Due to
the occurrence of integer part functions, it will be useful express
the splitting type $j$ of $E$ as $j = nk+b$, with $0 \leq b < k$.

\begin{corollary}\label{cor.loceuler}
Let $E$ be a rank-$2$ bundle over $Z_k$ of splitting type $j$ with $j
> 0$ and let $j=nk+b$ as above. The following are sharp bounds for the
local holomorphic Euler characteristic of $E$:
\[ j - 1 \leq \chi( \ell, E) \leq \begin{cases} n^2k+2nb+b-1 &
   \text{if } k \geq 2 \text{ and } 1 \leq b < k \text{ ,} \\
   n^2k & \text{if } k \geq 2 \text{ and } b = 0 \text{ ,} \end{cases} \]
and
\[ j \leq \chi(\ell, E) \leq j^2 \text{ for } k=1 \text{ .} \]
\end{corollary}
\begin{proof}
It follows directly from the definitions that for a curve $\ell^2=-k$
inside a smooth surface, $h^0\bigl(X;
(\pi_*E)^{\vee\vee}\bigl/\sigma_* E\bigr) = \bw_k(E)$ and $h^0(X;
R^1\sigma_*{E}) = \bh_k(E)$, so the result follows by adding the
results of Theorems \ref{thm.bndheight} and \ref{thm.bndwidth}.

To aid the computation, note that for $2 \leq b < k$ we have
$\bigl\lfloor \frac{j}{k} \bigr\rfloor = \bigl\lfloor \frac{j-2}{k}
\bigr\rfloor = n$, while for $b=0,1$ we have $\bigl\lfloor \frac{j}{k}
\bigr\rfloor = n$ and $\bigl\lfloor \frac{j-2}{k} \bigr\rfloor = n -
1$. We find that for $b=0,1$, $\chi( \ell, E) \leq n^2k+2nb$, and we
can absorb the case $b=1$ into the other case.

The upper bound is attained by the split bundle
$E=\Oh(-j)\oplus\Oh(j)$, while the lower bound is attained by the
``generic'' bundle with extension $p(z,u)=zu$ as in the proof of Lemma~\ref{lem.width_split}.
\end{proof}

We intend to study numerical invariants of bundles over more general
exceptional loci in future papers. One property of local holomorphic
Euler characteristics of \emph{instantons} (where $j=nk$ for some $n$)
is immediate from Corollary~\ref{cor.loceuler}:

\begin{corollary}
For a rank-$2$ bundle $E$ on $Z_k$ with $c_1(E)=0$ and splitting type
$j=nk$, the integer ranges $[1, \dotsc, k-2]$ and $[k+1, \dotsc,
2k-2]$ cannot occur as local holomorphic Euler characteristics
(``instanton charges''). These ranges are non-empty when $k>2$. (E.g.\
on $Z_3$, charges $1$ and $4$ cannot occur.\footnote{But the bundle
with $p(z,u)=z^{-1}u+z^4u^2$ has charge $7$. Also, it seems always
possible to create a new bundle with charge $+k$ via elementary
transformations, and since charges $2$ and $3$ exist by
Corollary~\ref{cor.loceuler}, we expect these are the only charge gaps
on $Z_3$.})
\end{corollary}

\section{Balancing}\label{sec.balancing}

We now consider the question of constructing vector bundles with
specified numerical invariants. We use the technique of
\emph{balancing bundles}. This technique was
used in \cite{FM} to prove the existence of bundles over $Z_1$ with
any prescribed numerically admissible invariants, and in \cite{RMJM}
to study bundles on $Z_2$.

Given two bundles $E$ and $E'$ of splitting type $(j_1, \dotsc, j_r)$
and $(j'_1, \dotsc, j'_r)$, we say that $E$ is \emph{more balanced}
than $E'$ if $j_1 -j_r \leq j'_1 - j'_r$. The advantage of balancing a
bundle is that we control the numerical invariants at each step, and
we only need to compute numerical invariants for a smaller range of
bundles.

The simplest case of balancing is for rank-$2$ bundles and goes as
follows. If $j_1 - j_2 \leq k-1$, we have won and we stop. If $j_2
\leq j_1 - k$ we make the construction of \cite{RMJM}, namely an
elementary transformation with respect to $\Oh_\ell(j_2)$ and obtain a
new more balanced bundle with splitting type $(j_1, j_2+k)$. We may
also compare the invariants of the two bundles; at the end we reduce
to a case with $j_1 - j_2 \leq k-1$. We now describe how to balance
bundles of rank $r \geq 2$. Let $E$ be a rank-$r$ vector bundle over
$Z_k$ of splitting type $j_1 \geq j_2 \geq \dotsb \geq j_r$. We say
that $E$ is \emph{balanced} if $j_1 \leq j_r + k - 1$. The objective
is to balance $E$. Balancing associates to $E$ the following data:
\begin{enumerate}
\item A positive integer $t$ (the number of steps);
\item a finite sequence of $r$-tuples of non-increasing integers
      $\{j(i,l)\}$ with $1 \leq i \leq t$ and $1 \leq l \leq r$ (the
      splitting types) satisfying:
      \begin{enumerate}
      \renewcommand{\labelenumi}{(\roman{enumi})}
      \item $j(1,l)=j_l$ for $1 \leq l \leq r$ (the splitting of the bundle $E$),
      \item $\sum_{l=1}^{r} j(i,l) = \sum_{i=1}^{r} j(1,l) + ki - k$
            for $2 \leq i \leq t$ (change of splitting produced by an
            elementary transformation), and
      \item $j(t,1) \leq j(t,r) + k - 1$ (arrive at a balanced bundle);
      \end{enumerate}
\item a chain $E_1, \dotsc, E_t$ of vector bundles, starting with
      $E_1 = E$, where $E_i$ has splitting type $j(i,1) \geq
      j(i,2) \geq \dotsb \geq j(i,r)$ for $1 \leq i \leq t$.
\end{enumerate}

\begin{defn}\label{dfn.admissible}
A sequence $\{j(i,l)\}$ of $r$-tuples of integers satisfying the
numerical properties 1, 2 and 3 above is called an \emph{admissible
sequence}. The sequence $\{j(i,l)\}_{i=1}^{t}$ of splitting types of
the bundles $E_i$ obtained in balancing the bundle $E$ is then called
the admissible sequence \emph{associated} to $E$.
\end{defn}

Balancing a bundle $E$ proceeds as follows: Set $E_1 \ce E$ and
$j(1,l) = j_l$ for $1 \leq l \leq r$. If $j(1,1) \leq j(1,r) + k - 1$
we have won, we set $t \ce 1$ and stop. Otherwise, $j(1,1) \geq j(1,r)
+ k$. Choose a surjective homomorphism $\rho \colon E\rest\ell \to
\Oh_\ell(j_r)$ and make the corresponding elementary transformation
$\br \colon E \to\Oh_\ell(j_r)$. Set $E_2 \ce \ker\br$. Since $j_r
\leq j_{r-1}$, we have $\ker\rho \iso \bigoplus_{i=1}^{r-1}
\Oh_\ell(j_i)$, and $\ker\br\rest\ell$ fits in the exact sequence
\begin{equation}\label{elm}
  0 \to \Oh_\ell(j_r+k) \to \ker\br\rest\ell \to \bigoplus_{i=1}^{r-1} \Oh_\ell(j_i) \to 0 \text{ .}
\end{equation}
We call $j(2,1) \geq j(2,2) \geq \dotsb \geq j(2,r)$ the
\emph{splitting type} of $E_2$. In particular,
\[ \sum_{l=1}^{r} j(2,l) = \sum_{l=1}^{r} j(1,l)+k \text{ .} \]
If $j(2,1) \leq j(2,r) + k - 1$, then we have won and we stop. If
$j(2,1) \geq j(2,r) + k$, then we obtain in the same way a bundle
$E_3$ with splitting type $j(3,1) \geq \dotsb \geq j(3,r)$ such that
$j(3,1) \leq j(2,1) \leq j(1,1)$, $j(3,r) \geq j(2,r) \geq j(1,r)$. We
continue this process, thus obtaining bundles $E_2, E_3, \dotsc E_i,
\dotsc$ with splitting type $j(i,1) \geq \dotsb \geq j(i,r)$, $j(i,1)
\leq j(i-1,1)$, $j(i,r) \geq j(i-1,r)$ and $\sum _{l=1}^{r} j(i,l) =
\sum _{l=1}^{r} j(1,l) + ki - k$ for all integers $i \geq 2$ for which
the bundle $E_i$ is defined. We must show that this procedure stops,
i.e.\ that after finitely many steps we arrive at a bundle $E_i$ such
that $j(i,1) \leq j(i,r) + k - 1$. At the same time we will show that,
in a suitable sense, each bundle $E_x$, $x \geq 2$, which is defined
is more balanced than the preceding one $E_{x-1}$. Call $a_1 \geq
\dotsb \geq a_r$ the splitting type of the bundle $T \ce
\mathcal{O}_\ell(j_r+k) \oplus \bigoplus_{i=1}^{r-1} \mathcal
{O}_\ell(j_i)$. Since $j_1 \geq j_r + k$ by assumption, we have $a_1 =
j_1$. If $j_{r-1} \geq j_r + k$, then $a_r = j_r + k$ and hence $a_1 -
a_r < j_1 - j_r$. If $j_{r-1} \leq j_r + k-1$, then $a_r = j_{r-1}$
and hence $a_1 - a_r = j_1 - j_{r-1}$. Hence in the latter case we
have $a_1 - a_r < j_1 - j_r$, unless $j_{r-1} = j_r$. By the exact
sequence \eqref{elm} and semi-continuity, the bundle $E_2\rest\ell$ is
more balanced (in the sense of the Harder--Narasimhan filtration) than
the bundle $T$. Hence we always have $j(2,1) - j(2,r) \leq j(1,1) -
j(1,r)$, while $j(2,1) - j(2,r) < j(1,1) - j(1,r)$ unless $j(1,r-1) =
j(1,r)$. Even in this case, since $\sum_{l=1}^{r} j(x,l) = k +
\sum_{l=1}^{r} j(x-1,l)$, we see that the procedure must stop after
finitely many steps.

We now construct families of bundles having prescribed associated
admissible sequence, generalising \cite[Theorem 2.1]{RMJM}, which
constructed such families for rank-$2$ bundles near a $-1$-curve. The
generalisation requires only minor modifications, but we give the
details for completeness.

\begin{theorem}\label{thm.admissible}
Fix an admissible sequence $\{j(i,l)\}_{i=1}^{t}$ and let $E$ and $F$
be rank-$r$ vector bundles on $\Chat$ with $\{j(i,l)\}_{i=1}^{t}$ as
an associated admissible sequence. Then there exists a flat family
$\{E_s\}_{s\in T}$ of rank-$r$ vector bundles on $\Chat$ parametrised
by an integral variety $T$ and $s_0, s_1 \in T$ with $E_{s_0}\iso E$
and $E_{s_1} \iso F$ such that $E_s$ has $\{j(i,l)\}_{i=1}^{t}$ as
admissible sequence for every $s \in T$.
\end{theorem}
\begin{proof}
We use induction on $t$. If $t=1$ the result is obvious because $E$
and $F$ are split vector bundles with the same splitting type and are
hence isomorphic. Assume that $t>1$ and that the result is true for
$t-1$. Let $E_2, F_2$ be the second bundles associated to $E, F$
respectively. Hence $E_2$ and $F_2$ have $\{j(2,l)\}_{i=2}^{t}$ as an
associated admissible sequence. By induction, there is a flat family
$\{E'_s\}_{s \in S}$ of rank-$r$ vector bundles on $\Chat$ and $m_0,
m_1 \in S$ with $E'_{m_0} \iso E_2$, $E'_{m_1} \iso F_2$ and such that
$E'_s$ has $\{j(i,l)\}_{i=2}^{t}$ as an associated admissible sequence
for every $s \in S$. We write $j'_i = j(2,i)$ to simplify notation.
By the balancing construction, the bundles $E$ and $F$ fit into exact
sequences
\begin{gather*}
  0 \to E \to E_2(\ell) \to \Oh_\ell(j_1-k) \to 0 \\
  0 \to F \to F_2(\ell) \to \Oh_\ell(j_1-k) \to 0
\end{gather*}
For every bundle $M$ on $\Chat$ having $\{j(i,l)\}_{i=2}^{t}$ as an
associated admissible sequence, the set of surjective homomorphisms
$\bt \colon M(\ell) \to \Oh_\ell(j_1-k)$ is parametrised by an
integral variety whose dimension depends only on $j_1, j_2$ and $j'_2
= j_1 + j_2 - j'_1 + k$. The kernel $\ker\bt\rest\ell$ is an extension
of $\Oh_\ell(j'_1)$ by $\Oh_\ell(j_1)$. This extension splits since
$j_1 \geq j'_1 + k$, and hence the bundle $\ker\bt$ has
$\{j(i,l)\}_{i=1}^{t}$ as an admissible sequence. Varying $M$ among
bundles $E'_s$ for $s \in S$ we get that the set of all such
surjections is parametrised by an irreducible non-empty variety
$T$. For any fixed ample line bundle $H$ on the $n^\mathrm{th}$
neighbourhood of $\ell$, it follows from the exact sequences in the
balancing construction that the bundles in this family have the same
Hilbert polynomial with respect to $H$, and therefore the family is
flat.
\end{proof}

\section{Moduli}\label{sec.moduli}

\begin{defn}
We write
\[ \sM_j(k) = \bigl\{E \to Z_k : E\rest\ell \iso \Oh(j) \oplus \Oh(-j)\bigr\} \bigl/ \sim \]
for the moduli stack of bundles over $Z_k$ of splitting type $j$.
\end{defn}

Recall from \eqref{eq!ext} that a bundle $E$ on $Z_k$ of splitting
type $j$ is an extension of $\Oh(j)$ by $\Oh(-j)$ and is therefore
determined by its extension class. In our choice of coordinates, this
amounts to saying that $E$ is determined by the pair $(j,p)$, where
$j$ is the splitting type and $p$ is a polynomial as in
\eqref{poly}. Let $N \ce \bigl\lfloor \frac{2j-2}k \bigr\rfloor$ as
before; then $p$ has $m \ce N(2j-1) - k \binom{N+1}{2}$
coefficients. We identify $p$ as an element in $\C^m$ by writing its
coefficients in lexicographical order. We then define the equivalence
relation $p \sim p'$ if $(j,p)$ and $(j,p')$ define isomorphic bundles
over $Z_k$. We give $\C^m$ the quotient topology. There is a bijection
\begin{eqnarray*}
  \phi \colon \sM_j(k) &\to& \C^m\bigl/\sim \text{ ,} \\
   \begin{pmatrix} z^{j} & p \\ 0 & z^{-j}\end{pmatrix}
   &\mapsto& \bigl(p_{1,k-j+1}, \dotsc, p_{N,j-1}\bigr) \text{ ,}
\end{eqnarray*}
where 
\[ p(z,u) = \sum_{r=1\vphantom{jks}}^N \ \  \sum_{s = kr -j + 1}^{j - 1}
   p_{rs} \, z^s u^r \text{ .} \]
We give $\sM_j(k)$ the topology induced by this bijection. Here are
some examples:

\begin{exa}
For each $k$, $\sM_0(k)$ contains only one point, corresponding to the
trivial bundle over $Z_k$. In other words, if a bundle over $Z_k$ is
trivial over the zero section, it is globally trivial.
\end{exa}

\begin{exa}
For each $2j-2<k$, $\sM_j(k)$ contains only one point. In other words,
a holomorphic bundle over $Z_k$ of splitting type $2j-2<k$
splits. This can be verified directly from formula \eqref{poly}.
\end{exa}

\begin{exa}
$\sM_2(2)$ contains exactly two points, \cite[Theorem 6.24]{SO}.
\end{exa}

\begin{exa}
$\sM_2(1) \simeq \PP^1 \cup \{A,B\}$, where $A$ and $B$ are points,
with open sets $U \subset \PP^1$, where $U$ is open in the usual
topology of $\PP^1$, $\PP^1\cup\{A\}$, and the whole space,
\cite[Theorem~4.2]{JA}.
\end{exa}

\begin{exa}
$\sM_3(2) \simeq \PP^2 \cup \{A,B\}$, where $A$ and $B$ are points,
with open sets $U \subset \PP^2$, where $U$ is open in the usual
topology of $\PP^2$, $\PP^2 \cup \{A\}$, and the whole space,
\cite[Theorem~6.35]{SO}.
\end{exa}

\begin{exa}
$\sM_j(k)$ is non-Hausdorff for $2j-2 \geq k$. This uses Theorem
\ref{thm.genopset} below to show non-emptiness (positive
dimensionality if the inequality is strict). Notice that the only open
neighbourhood of the split bundle is the entire $\sM_j(k)$.
\end{exa}

\begin{lemma}\label{lem.proj}
If $p' = \lambda p$ for some $\lambda \in \mathbb{C}^\times$, then the matrices
\[ \begin{pmatrix} z^j & p \\ 0 & z^{-j} \end{pmatrix} \quad\text{and}\quad
   \begin{pmatrix} z^j & p'\\ 0 & z^{-j} \end{pmatrix} \]
give holomorphically equivalent vector bundles.
\end{lemma}
\begin{proof}
Just write down the isomorphism as
\[ \begin{pmatrix} z^j & p' \\ 0 & z^{-j} \end{pmatrix} = 
   \begin{pmatrix} 1 & 0 \\ 0 & 1/\lambda \end{pmatrix}
   \begin{pmatrix} z^j & p \\ 0 & z^{-j} \end{pmatrix}
   \begin{pmatrix} 1 & 0 \\ 0 & \lambda \end{pmatrix} \text{ .} \]
\end{proof}

\paragraph{Notation:} We write $\mathcal{F}^{(n)} \ce
\mathcal{F}\rvert_{\ell_n} = \mathcal{F} \otimes
\mathcal{O}\bigl/\mathcal{I}_\ell^{n+1}$ for any coherent sheaf
$\mathcal{F}$ and similarly $E^{(n)} \ce E\rvert_{\ell_n}$ for vector
bundles on $Z_k$.

\begin{theorem}\label{thm.proj}
On the first infinitesimal neighbourhood, two bundles $E^{(1)}$ and
$E^{\prime(1)}$ with respective transition matrices
\[ \begin{pmatrix} z^j & p_1 \\ 0 & z^{-j} \end{pmatrix} \quad\text{and}\quad
   \begin{pmatrix} z^j & p'_1 \\ 0 & z^{-j} \end{pmatrix} \]
are isomorphic if and only if $p^\prime_1 = \lambda p_1$ for some
$\lambda \in \mathbb{C}^\times$.
\end{theorem}
\begin{proof}
The ``if'' part is Lemma \ref{lem.proj}. Now suppose $E^{(1)}$ and
$E^{\prime(1)}$ are isomorphic. According to our notation we have $p_1
= \sum_{s = k-j+1}^{j-1} p_{1s}\,z^s\,u$ and $p'_1 = \sum_{s =
k-j+1}^{j-1} p'_{1s}\,z^l\,u$. We will write the isomorphism in the
form
\[ \begin{pmatrix}z^j & p'_1 \\ 0 & z^{-j} \end{pmatrix}
   \begin{pmatrix} a & b \\ c & d \end{pmatrix} = 
   \begin{pmatrix} \alpha & \beta \\ \gamma & \delta \end{pmatrix}
   \begin{pmatrix} z^j & p_1 \\ 0 & z^{-j} \end{pmatrix} \text{ ,} \]
where $a$, $b$, $c$ and $d$ are holomorphic on $U$, and $\alpha$,
$\beta$, $\gamma$ and $\delta$ are holomorphic on $V$. On the first
infinitesimal neighbourhood, this yields the following set of equations:
\[ \left\{ \begin{array}{rclr}
   \bigl(a_0(z) + a_1(z)u \bigr) z^j + p'_1 c_0(z)
   & = & \bigl(\alpha_0(z^{-1}) + \alpha_1(z^{-1})zu\bigr) z^j & \qquad(\text{A1}) \\
   z^{-j} \bigl(c_0(z) + c_1(z)u\bigr)
   & = & \bigl(\gamma_0(z^{-1}) + \gamma_1(z^{-1})zu\bigr) z^j & (\text{A2}) \\
   \bigl( b_0(z) + b_1(z)u \bigr) z^j + p'_1 d_0(z)
   & = & \alpha_0(z^{-1}) p_1 + \bigl(\beta_0(z^{-1}) + \beta_1(z^{-1})zu\bigr) z^{-j} & (\text{A3}) \\
   z^{-j} \bigl(d_0(z) + d_1(z)u \bigr) 
   & = & \gamma_0(z^{-1}) p_1 + \bigl(\delta_0(z^{-1}) + \delta_ 1(z^{-1})zu\bigr) z^{-j} & (\text{A4})
   \end{array} \right. \]
Recalling that $p_1$ and $p^\prime_1$ are multiples of $u$ and
equating terms that are independent of $u$ in (A1) and (A4) gives
$a_0(z) = \alpha_0(z^{-1})$ and $d_0(z) = \delta_0(z^{-1})$
respectively. Therefore $a_0$, $\alpha_0$, $d_0$ and $\delta_0$ are
constants, and $a_0 = \alpha_0$ and $d_0 = \delta_0$. Next we equate
terms in $u$ in Equation (A3), obtaining
\[ b_1(z)u \, z^j + p'_1 d_ 0 = \alpha_0  p_1 + \beta_1(z^{-1}) u \, z^{-j} \text{ .} \]
In Equation (A3), $z^j\,b_1$ has only terms $z^l$ for $l \geq j$, and
$z^{-j}\,\beta_1$ has only terms $z^l$ for $l \leq -j$. Consequently,
they do not affect the terms appearing in $p_1$ and $p'_1$; and the
remaining part of Equation (A3) gives $p'_1 d_0 = \alpha_0 p_1$.  We
observe that $p_1$ and $p'_1$ differ by a constant.

It remains to show that $d_0$ and $\alpha_0$ are non-zero. Taking
terms that are independent of $u$ in Equation (A3) we have
$b_0(z)\,z^j = \beta_0(z^{-j})\,z^{-j}$, which implies $b_0(z) =
\beta_0(z^{-1}) = 0$. It follows that over the exceptional divisor our
coordinate change has determinant $a_0\, d_0$, hence
$\alpha_0\,\delta_0 = a_0\,d_0 \neq 0$.
\end{proof}

\begin{rem}[Construction of the generic set and further strata]
For a fixed bundle $E$ over $Z_k$ with transition matrix
$\bmat{z^j}{p}{0}{z^{-j}}$, we write any automorphism of $E$ in the
form
\begin{equation}\label{eq.aut}
  \begin{pmatrix} \alpha & \beta \\ \gamma & \delta \end{pmatrix}
  = \begin{pmatrix} z^{j} & p \\ 0 & z^{-j} \end{pmatrix}
    \begin{pmatrix} a & b \\ c & d \end{pmatrix}
    \begin{pmatrix} z^{-j} & -p \\ 0 & z^{j} \end{pmatrix}
  \text{ ,}
\end{equation}
where $a$, $b$, $c$, $d$ are holomorphic on $U$ and $\alpha$, $\beta$,
$\gamma$, $\delta$ are holomorphic on $V$. This gives rise to a set of
four equations, imposing conditions on $a,b,c,d$ that make the
expression on the right-hand side holomorphic on $V$ (see the Long
Proof of Theorem \ref{thm.genopset}). Therefore, the set of
automorphisms of $E$ is determined by the matrix $\mat{a}{b}{c}{d}$
subjected to certain constraints for holomorphicity. The number of
constraints depends on $p$, and for a fixed $j$ there is a minimal
number of constraints that happens for ``most'' of the bundles. We
call {\em generic} a bundle whose automorphisms realise the minimal
number of constraints. It would be desirable to phrase this as ``a
bundle is generic (or, even better, stable) when it has the smallest
number of automorphisms'' as in the compact case; alas, the spaces of
automorphism for each bundle is infinite-dimensional; moreover, no
concept of stability is available for bundles on non-compact spaces.

Take for example the extreme cases: for a fixed $j$ we have $p(z,u) =
\sum p_{rs} \, z^s u^r$, and the two opposite cases are (i) when
$p_{rs} \neq 0 \ \forall r,s$ and the corresponding bundle is generic;
and (ii) when $p=0$ and the corresponding bundle belongs to the least
generic stratum. What is happening is the simple fact that extra
conditions on $p$ make the equations on $a,b,c,d$ harder to solve.

One can see the division of $\mathcal M_j(k)$ into strata as a
constructive algorithm that proceeds by infinitesimal
neighbourhoods. Start with the first infinitesimal neighbourhood,
where we know that the only automorphisms are scalar
multiplication. Now, pass to the second infinitesimal neighbourhood
and keep in the generic stratum all bundles having the most
automorphisms, that is, the ones whose expression of $p$ is such that
the minimal number of constraints is imposed on $a,b,c,d.$ Separate
the other bundles away form the generic stratum, again divided by
number of automorphisms. Now restart the process on the third
neighbourhood, and so on. When this process is carried out on the
neighbourhood $N = \left\lfloor (2j - 2)/k \right\rfloor$, then we are
done, so it is a finite procedure. This process fixes the bundle
$\End(E)$ and therefore its cohomologies $H^i\bigl(\SEnd(E)\bigr)$. In
particular $H^1\bigl(\SEnd(E)\bigr)$ is fixed, and so is $H^1(E)$ and
correspondingly the height $\bh_k(E)$. Therefore only bundles with the
same height belong to the same stratum (but not conversely). The
generic stratum is contained in the set of bundles with minimal
$h^1\bigl(\SEnd(E)\bigr)$ (like stable points in the compact case),
and hence minimal height. Now, $H^0\bigl(\SEnd(E)\bigr)$ is also fixed
by the $N^\text{th}$ neighbourhood, and so is $H^0(E)$. This is not
quite the same as fixing the width, though. We have
\[ 0 \rightarrow S^n(N_{\ell,Z_k}^*) \rightarrow \Oh^{(n)} \rightarrow
   \Oh^{(n-1)} \rightarrow 0 \text{ ,} \]
where we write $N_{\ell,Z_k}$ for the normal sheaf of $\ell \subset
Z_k$, and tensoring with $E$ gives
\[ 0 \rightarrow S^n(N_{\ell,Z_k}^*)\otimes E \rightarrow E^{(n)}
   \rightarrow E^{(n-1)} \rightarrow 0 \text{ .} \]
Now taking the direct image under $\pi$ gives
\[ 0 \rightarrow \pi_*\bigl(S^n(N_{\ell,Z_k}^*)\otimes E\bigr) \rightarrow 
   \pi_*\bigl(E^{(n)}\bigr) \rightarrow \pi_*\bigl(E^{(n-1)}\bigr) \rightarrow
   R^1\pi_* \bigl(S^n(N_{\ell,Z_k}^*)\otimes E\bigr) \rightarrow \dotsb \text{ .} \]
Now for sufficiently high $n$, $\pi_*E^{(n)}$ completely determines
the width. We can compare the corresponding matrices of endomorphisms
for two bundles $E$ and $E'$. Suppose they have the same width, but yet
$H^0\bigl(\SEnd(E)^{(n)}\bigr)$ and $H^0\bigl(\SEnd(E')^{(n)}\bigr)$
have positive relative dimension. Here, because the canonical
expressions of extension classes are cut off at $u$-degree $N =
\left\lfloor (2j - 2)/k \right\rfloor$, we cannot observe differences
happening at intervals smaller than $k$, so we must now restrict
ourselves to the instanton case, where $j=nk$, and correspondingly we
can observe the relative values of $h^0$ properly. Comparing the above
sequence with $\SEnd(E)$ and $\SEnd(E')$ in place of $E$, we see that
the only way this can happen is that also the $R^1\pi_*$-terms have
different relative dimensions, but in this case the heights are
different, and the bundles belong to different strata.
\end{rem}

\begin{theorem}\label{thm.genopset}
For $j \geq k$, $\sM_j(k)$ has an open, dense subspace homeomorphic
to a complex projective space \ $\PP^{2j-2-k}$ minus a closed
subvariety of codimension at least $k+1$.
\end{theorem}
\begin{proof}[Short proof]
\cite[Theorem 3.5]{JA} showed that the generic set of $\sM_j(1)$ is a
projective space $\PP^{2j-3}$ minus a closed subvariety of codimension
$\geq 2$. The only modification needed to generalise the proof to
$k>1$ is the calculation of dimension of the generic set. Generic
bundles do not split on the first infinitesimal neighbourhood, and
there the only equivalence relation is projectivisation, by Theorem
\ref{thm.proj}. The dimension count follows from formula \eqref{poly},
which shows that the $u$-coefficients are $\sum_{s=k-j+1}^{j-1}
p_{1s}$. There are $2j-k-1$ coefficients, and after projectivising we
obtain $\PP^{2j-k-2}$. However, not all points of this $\PP^{2j-k-2}$
are generic, and one must remove a closed subvariety where too much
vanishing of coefficients occurs, thus implying larger numerical
invariants. The closed subvariety to be removed contains all points
having coefficients $z^s u = 0$ for $0 < s < k$, hence is defined by
at least $k+1$ equations.
\end{proof}
\begin{proof}[Long proof]
Let $E$ and $E'$ be the bundles on $Z_k$ given respectively by the
extension classes $p$ and $p'$. By Theorem \ref{thm.proj}
we know that on $\ell_1$ the only isomorphism of
bundles is scaling, so we may assume
\[ p = p_1 + p_2 \qquad\text{and}\qquad p' = p_1 + p'_2 \text{ ,} \]
where $p_1 \ce p\rvert_{\ell_1}$. If $E$ and $E'$ are isomorphic, then
\begin{multline}\label{eq.iso}
  \begin{pmatrix} \alpha & \beta \\ \gamma & \delta \end{pmatrix}
  = \begin{pmatrix} z^{j} & p' \\ 0 & z^{-j} \end{pmatrix}
    \begin{pmatrix} a & b \\ c & d \end{pmatrix}
    \begin{pmatrix} z^{-j} & -p \\ 0 & z^{j} \end{pmatrix} 
  = \begin{pmatrix} a + z^{-j} p' c & z^{2j}b + z^j(dp'-ap) - cpp' \\
  z^{-2j} c & d - z^{-j}pc \end{pmatrix} \text{ ,}
\end{multline}
for some change of coordinates $\mat{a}{b}{c}{d}$ holomorphic on $U$
and $\mat{\alpha}{\beta}{\gamma}{\delta}$ holomorphic on $V$, which we
may assume have determinant one. We write
\[ \alpha(z^{-1},z^ku) = \alpha_0(z^{-1}) + \alpha_1(z^{-1}) \, z^ku +
   \dotsb \text{ , similarly for $\beta$, $\gamma$, $\delta$, and} \]
\[ a(z,u) = a_0(z) + a_1(z)u + a_2(z)u^2 + \dotsb \text{ , similarly for $b$, $c$, $d$,} \]
where the coefficients are convergent power series in $z^{-1}$ or $z$,
respectively.
Then
\[ \begin{pmatrix} \alpha_0 & \beta_0 \\ \gamma_0 & \delta_0 \end{pmatrix}
   = \begin{pmatrix} a_0 & b_0 z^{2j} \\ c_0 z^{-2j} & d_0 \end{pmatrix}
  \text{ ,} \]
implies $\beta_0 = b_0 = 0$, and $c_0(z) = c_{00} + c_{01}z + \dotsb +
c_{0,2j}z^{2j}$. It follows that $\alpha_{00} = a_{00} = \lambda =
\delta_{00}^{-1} = d_{00}^{-1} = 1$.

The coefficients of $u$ are
\begin{equation}\label{eq.isol1}
  \begin{pmatrix} \alpha_1 z^ku & \beta_1 z^ku \\[1ex] \gamma_1 z^ku & \delta_1 z^ku \end{pmatrix}
  = \begin{pmatrix} a_1 u + p_1 c_0z^{-j} & b_1 u z^{2j} + z^j\cancel{(d_0p_1-a_0p_1)}\\[1ex]
  c_1 u z^{-2j} & d_1 u - p_1c_0z^{-j} \end{pmatrix} \text{ ,}
\end{equation}
which has to be holomorphic in $(z^{-1},z^ku)$. This forces $b_1(z) =
b_{10} + b_{11}z + \dotsb + b_{1,-2j+k}z^{-2j+k}$, whence $b_1 \neq 0$
provided $k - 2j \geq 0$, i.e.\ $j \leq \lfloor \frac{k}{2}
\rfloor$. But by assumption, $j \geq k$, so that $b_1 = 0$, and
consequently $d_1 = -a_1$. Furthermore, assuming for the moment that
$j \geq k+1$, the $(1,1)$-entry of \eqref{eq.isol1} entails the
following relations between the terms of $a_1$ and $c_0$:
\begin{equation}\label{eq.acrel}
  \begin{pmatrix} a_{1,k+1}\\a_{1,k+2}\\ \vdots \\ a_{1,2j-2} \\ a_{1,2j-1} \end{pmatrix}
  + \begin{pmatrix} p_{1,j-1} & p_{1,j-2} & \dotso & p_{1,k-2+2} & p_{k-j+1} \\ 0 & p_{1,j-1} &
  \dotso & p_{1,k-j+3} & p_{1,k-j+2} \\ \vdots & \vdots & \ddots & \vdots & \vdots \\
  0 & 0 & \dotso & p_{1,j-1} & p_{1,j-2} \\ 0 & 0 & \dotso & 0 & p_{1,j-1} \end{pmatrix}
  \begin{pmatrix} c_{0,k+2} \\ c_{0,k+3} \\ \vdots \\ c_{0,2j-1} \\ c_{0,2j} \end{pmatrix}
  = 0 \text{ ,}
\end{equation}
and $a_{1,s}=0$ for $s\geq 2j$.

On $\ell_2$, the terms in \eqref{eq.iso} with $u^2$ are 
\begin{multline*}
  \begin{pmatrix} \alpha_2 z^{2k} u^2 & \beta_2 z^{2k} u^2 \\[1ex] \gamma_2 z^{2k} u^2 & \delta_2 z^{2k} u^2 \end{pmatrix}
  = \\ 
  \begin{pmatrix} a_2 u^2 + z^{-j}(p'_2\,c_0 + p_1\,c_1\,u) &
  b_2\,u^2\,z^{2j} + z^j\bigl( (d_1 - a_1) \, u \, p_1 + (p_2-p'_2) \bigr) - c_0\,p_1^2 \\[1ex]
  c_2\,u^2\,z^{-2j} & d_2\,u^2 - z^{-j}(p_2\,c_0 + p_1\,c_1\,u) \end{pmatrix} \text{ .}
\end{multline*}
We need to  examine the $(1,2)$-entry:
\begin{equation}\label{eq.l2-12}
  b_2\,u^2\,z^{2j} + z^j\bigl( (d_1 - a_1) \, u \, p_1 + (p_2-p'_2) \bigr) - c_0\,p_1^2
\end{equation}

The conditions on the expression \eqref{eq.l2-12} is that all
coefficients of $z^{l}\,u^2$ vanish for $l \geq 2k$, but any coefficient
of $z^{2j}\,u^2$ and higher can be cancelled by choosing $b_2$
appropriately. So we only need to consider the range $k+1 \leq l
\leq 2j-1$ to verify when the expression
\[  z^j (d_1-a_1) u\,p_1 + z^j(p\low_2-p'_2) - c_0\,p_1^2  \]
can be made holomorphic on $V$ for any choice of $p'_2$. 
Given the determinant-one condition on the coordinate changes, this becomes
\[ z^j (-2a_1) u\,p_1 + z^j(p\low_2-p'_2) - c_0\,p_1^2 \text{ ,} \]
and plugging in the values of $p$,
\begin{equation}
  -2 z^j a_1 \sum_{s=k-j+1}^{j-1} p\low_{1s}\,z^s\,u +
  z^j \sum_{s=2k-j+1}^{j-1}\bigl(p\low_{2,s}-p'_{2,s}\bigr)z^s\,u^2 - c_0
  \left( \sum_{s=k-j+1}^{j-1} p\low_{1s}\,z^s\,u\right)^2 \text{ .}
\end{equation}
The terms to be  cancelled are:
\begin{itemize}
\item Step $1$, the coefficient of $z^{2k+1}u^2$: 
      \begin{multline*}
        -2 \bigl(a\low_{1,k} p\low_{1,k-j+1} + a\low_{1,k-1} p\low_{1,k-j+2} + \dotsb +
        a\low_{1,0} p\low_{1,2k-j+1}\bigr) + \bigl(p\low_{2,2k-j+1} - p'_{2,2k-j+1}\bigr) \\
        + c\low_{0,2j-1} \, p^2_{1,k-j+1} + 2 c\low_{0,2j-2} \, p\low_{1,k-j+1} \, p\low_{1,k-j+2} +
        c\low_{0,2j-3} \bigl(p^2_{1,k-j+2} + 2 p\low_{1,k-j+1} \, p\low_{1,k-j+3}\bigr) \\
        + \dotsb + c\low_{0,1} \bigl(p^2_{1,k} + 2 p\low_{1,k-1} \, p\low_{1,k+1} + \dotsb \bigr)
        + c\low_{0,0} \bigl(2 p\low_{1,k} \, p\low_{1,k+1} + 2 p\low_{1,k-1} \, p\low_{1,k+2} + \dotsb\bigr) = 0
      \end{multline*}
\item Step $2$, the coefficient of $z^{2k+2}u^2$: 
      \begin{multline*}
        -2 \bigl(a\low_{1,k+1} p\low_{1,k-j+1} + a\low_{1,k} p\low_{1,k-j+2} + \dotsb +
        a\low_{1,0} p\low_{1,2k-j+2}\bigr) + \bigl(p\low_{2,2k-j+2} - p'_{2,2k-j+2}\bigr) \\
        + c\low_{0,2j} \, p^2_{1,k-j+1} + 2 c\low_{0,2j-1} \, p\low_{1,k-j+1} \, p\low_{1,k-j+2} +
        c\low_{0,2j-2} \bigl(p^2_{1,k-j+2} + 2 p\low_{1,k-j+1} \, p\low_{1,k-j+3}\bigr) \\
        + \dotsb + c\low_{0,0} \bigl(p_{1,k+1}^2 + 2 p\low_{1,k} \, p\low_{1,k+2} + \dotsb\bigr) = 0
      \end{multline*}
\item Step $s$, the coefficient of $z^{2k+s}u^2$ for $1 \leq s \leq 2j-2k-1$, until\ldots
\item Step $2j-2k-1$, the coefficient of $z^{2j-1}u^2$: 
      \begin{multline*}
        -2 \bigl(a\low_{1,2j-k-2} p\low_{1,k-j+1} + a\low_{1,2j-k-3} p\low_{1,k-j+2} + \dotsb +
        a\low_{1,0} p\low_{1,j-1}\bigr) + \bigl(p\low_{2,j-1} - p'_{2,j-1}\bigr) \\
        + c\low_{0,2j}\bigl(2 p\low_{1,-1} \, p\low_{1,0} + \dotsb \bigr) + c\low_{0,2j-1} \bigl(p^2_{1,0}
        + 2 p\low_{1,-1} \, p\low_{1,1} + \dotsb \bigr) + \dotsb + 
        c\low_{0,1} p^2_{1,j-1} = 0
      \end{multline*}
\end{itemize}
Now \emph{assume $k \leq j-1$}, that is, assume that $j$ is
large. Then a term in $p^2_{1,0}$ appears with some $c_{0,s}$ in each
of the above equations. Thus, choosing $c_0$ appropriately, we can
solve them all, and we conclude that there are only restrictions on
$p'_2$ when $p_{1,0} = 0$. Now, we can carry out a similar argument
for the coefficients $p_{1,s}$ for each $0 \leq s \leq k$, so the set
of non-generic bundles lives on the subvariety singled out by the
equations
\[ p_{1,0} = p_{1,1} = \dotsb = p_{1,k} = 0 \text{ ,} \]
thus having codimension at least $k+1$.

In the remaining case where $k=j$, we see directly from \eqref{poly}
that $p(z,u) = \sum_{s=1}^{k-1} p_{1s}\,z^s\,u$. Thus the only
non-generic bundle is the split bundle, and the generic set of
$\sM_k(k)$ is precisely $\PP^{k-2}$ when $k\geq2$, and empty when
$k=1$. (The same argument shows that the generic set is all of
$\PP^{2j-k-2}$ for $j \leq k \leq 2j-2$.)
\end{proof}

\begin{theorem}\label{embed}
There is a topological embedding \ $\Phi \colon \sM_j(k) \to
\sM_{j+k}(k)$. The image of \ $\Phi$ consists of all bundles in
$\sM_{j+k}(k)$ that split on the second infinitesimal neighbourhood of
$\ell$.
\end{theorem}
\begin{proof}
Using the identification $\phi \colon \sM_j(k) \to \C^m\bigl/\sim$,
we define a map
\[ \Phi \colon \sM_j(k) \to \sM_{j+k}(k) \quad\text{by}\quad
   (j,p) \mapsto(j+k, z^k u^2 p) \text{ .} \]
We want to show that $\Phi$ defines an embedding. We first show that the
map is well defined: Suppose that $\bmat{z^j}{p}{0}{z^{-j}}$ and
$\bmat{z^j}{p'}{0}{z^{-j}}$ represent isomorphic bundles. Then there are
coordinate changes $\bmat abcd$ holomorphic in $(z,u)$ and $\bmat\al\be\ga\de$
holomorphic in $(z^{-1}, z^k u)$ such that
\[ \begin{pmatrix}\al & \be \\ \ga & \de \end{pmatrix}
   = \begin{pmatrix}z^j & p' \\ 0 & z^{-j} \end{pmatrix}
     \begin{pmatrix}a & b \\ c & d \end{pmatrix}
     \begin{pmatrix}z^{-j} &-p \\ 0 & z^j \end{pmatrix} \text{ .} \]
Therefore these two bundles are isomorphic exactly when the system of
equations
\begin{equation}\label{eq.pr1}
  \begin{pmatrix}\al & \be \\ \ga & \de \end{pmatrix}
  = \begin{pmatrix}a+z^{-j}p'c & z^{2j}b+z^j(p'd-ap)-pp'c \\
  z^{-2j}c & d-z^{-j}pc \end{pmatrix}
\end{equation}
can be solved by a matrix $\mat abcd$ holomorphic in $(z,u)$ which makes
$\mat\al\be\ga\de$ holomorphic in $(z^{-1}, z^k u)$.

On the other hand, the images of these two bundles are given by
transition matrices $\mat{z^{j+k}}{z^k u^2 p}{0}{z^{-j-k}}$ and
$\mat{z^{j+k}}{z^k u^2 p'}{0}{z^{-j-k}}$, which represent isomorphic
bundles if and only if there are coordinate changes $\mat{\bar a}{\bar
b}{\bar c}{\bar d}$ holomorphic in $(z,u)$ and
$\mat{\bar\al}{\bar\be}{\bar\ga}{\bar\de}$ holomorphic in $(z^{-1},
z^k u)$ satisfying the equality
\[ \begin{pmatrix}\bar\al & \bar\be \\ \bar\ga & \bar\de \end{pmatrix}
   = \begin{pmatrix}z^{j+k} & z^ku^2p' \\ 0 & z^{-j-k} \end{pmatrix}
   \begin{pmatrix}\bar a & \bar b \\ \bar c & \bar d \end{pmatrix}
   \begin{pmatrix}z^{-j-k} &-z^ku^2p \\ 0 & z^{j+k} \end{pmatrix} \text{ .} \]
That is, the images represent isomorphic bundles if the system
\begin{equation}\label{eq.pr2}
   \begin{pmatrix}\bar\al & \bar\be \\ \bar\ga & \bar\de \end{pmatrix}
   = \begin{pmatrix}\bar a+z^{-j } u^2p'\bar c &
   z^{2k}\bigl(z^{2j}\bar b+z^{j}u^2(p'\bar d-\bar a p)-u^4pp'\bar c\bigr) \\
   z^{-2j-2k}\bar c & \bar d-z^{-j} u^2p\bar c \end{pmatrix}
\end{equation}
has a solution.

Write $x = \sum x_i u^i$ for $x \in \{a,b,c,d, \bar a,\bar b,\bar
c,\bar d\}$ and choose
\[ \bar a_i = a_{i+2k} \text{ ,}\quad \bar b_i = b_{i+2k}u^2 \text{ ,}\quad
   \bar c_i = c_{i+2k}u^{-2} \text{ ,}\quad \bar d_i = d_{i+2k} \text{ .} \]
Then if $\mat abcd$ solves \eqref{eq.pr1}, one verifies that $\mat{\bar a}{\bar
b}{\bar c}{\bar d}$ solves \eqref{eq.pr2}, which implies that the images
represent isomorphic bundles and therefore $\phi$ is well defined. To
show that the map is injective, just reverse the previous
argument. Continuity is obvious. Now we observe also that the image
$\phi({\sM}_j)$ is a saturated set in $\sM_{j+k}$, meaning that if
$y\sim x$ and $x\in\phi({\sM}_j)$ then $y\in\phi({\sM}_j)$. In fact,
if $E \in \phi({\sM}_j)$, then $E$ splits on the second infinitesimal
neighbourhood. Now if $E' \sim E$, then $E'$ must also split on the
second infinitesimal neighbourhood, and therefore the polynomial
corresponding to $E'$ is of the form $u^2p'$, and hence
$\phi(z^{-k}p')$ gives $E'$. Note also that $\phi({\sM}_j)$ is a
closed subset of ${\sM}_{j+k}$, given by the equations $p_{il} = 0$
for $i = 1, 2$. Now the fact that $\phi$ is a homeomorphism over its
image follows from the easy Lemma \ref{lem.homeo} given below.
\end{proof}

\begin{rem}
R.\ Moraru gave us the following coordinate-free expression of the
embedding map $\Phi \colon \sM_j(k) \to \sM_{j+k}(k)$:
\[ \Phi\bigl(E\bigr) = \otimes\,\Oh(-k) \circ \Elm_{\Oh_\ell(j+k)} \circ
   \Elm_{\Oh_\ell(j)} \bigl(E\bigr) \text{ ,} \]
where $\Elm_L$ denotes the elementary transformation with respect to
the line bundle $L$. Using this coordinate-free expression it becomes
obvious that $\Phi$ is well defined.
\end{rem}

\begin{lemma}\label{lem.homeo}
Let $X \subset Y$ be a closed subset and $\sim$ an equivalence
relation in $Y$ such that $X$ is $\sim$-saturated. Then the map $I
\colon X\bigl/\sim\to Y\bigl/\sim$ induced by the inclusion is a
homeomorphism over the image.
\end{lemma}
\begin{proof}
Denote the projections by $\pi_X \colon X \to X\bigl/\sim$ and $\pi_Y
\colon Y \to Y\bigl/\sim$. Let $F$ be a closed subset of
$X\bigl/\sim$. Then $\pi_X^{-1}(F)$ is closed and saturated in $X$,
and therefore $\pi_X^{-1}(F)$ is also closed and saturated in $Y$. It
follows that $\pi_Y\bigl( \pi_X^{-1}(F)\bigr)$ is closed in
$Y\bigl/\sim$.
\end{proof}

\begin{theorem}\label{thm.hdcomps}
If $j=nk$ for some $n \in \mathbb{N}$, then the pair $(\bh_k, \bw_k)$
stratifies instanton moduli stacks $\sM_j(k)$ into Hausdorff
components.
\end{theorem}
\begin{proof}
This proof uses the same techniques as that of \cite[Theorem
4.1]{PAMS}. On the first infinitesimal neighbourhood, we have two
possibilities: in the first case we have bundles belonging to the open
dense subset $\PP^{2j-2-k}$, singled out by having the lowest possible
values of numerical invariants (see Remark \ref{rem.stable} for a question
about stability); in the second case, at least one of the invariants
is strictly higher than the lower bound, and such bundles are
separated away from the most generic stratum. On the second
infinitesimal neighbourhood, the problem is solved by first separating
the most generic stratum from the other ones. For the remaining part
of the second neighbourhood, one divides the polynomial by $u$ falling
back to the same analysis done for the first neighbourhood. We are
then left only with bundles which split on the second
neighbourhood. We use induction $j$, assuming that the invariants
stratify $\sM_{j-k}(k)$ into Hausdorff components together with the
embedding Theorem \ref{embed}, stating that
$\Phi\bigl(\sM_{j-k}(k)\bigr)$ is the set of bundles on $\sM_{j}(k)$
that split on the second infinitesimal neighbourhood.
\end{proof}

\begin{rem}\label{rem.stable}
When this paper was nearly completed, we noticed the need to define a
notion of stability on $\sM_k(j)$ in order to have full compatibility
between the methods used in the proof of Theorem \ref{thm.genopset}
and the stratification presented in Theorem \ref{thm.hdcomps}. Defining
stability for a bundle $E$ via counting the dimension of $H^1(Z_k;\;\SEnd
E)$ seemed to us to be the most natural choice, and had the pleasant
feature of being a notion that fits extremely well with the standard
deformation theory that is well known for the compact case;
nevertheless, it was just an ad-hoc definition of stability. Since
there are many inequivalent ways to define stability, we chose to
postpone this question to a future paper.
\end{rem}

\begin{exa}\label{ex.nonhaus}
$\sM_3(1)$ is the simplest example in which the local holomorphic
Euler characteristic $\chi(\ell,E)$ does not distinguish Hausdorff
components: Writing $\chi = \bw_1 + \bh_1$, the (Hausdorff) generic
set has bundles with $\chi=3=1+2$, but the set of bundles with
$\chi=5$ is non-Hausdorff, containing both bundles of characteristic
$\chi=3+2$ and $\chi=2+3$. It is necessary to fix both $\bw_1$ and $\bh_1$
to obtain Hausdorff subspaces of $\sM_3(1)$.
\end{exa}

\begin{exa}\label{ex.noninst}
Bundles on $Z_2$ with splitting type $j=3$ do not represent
instantons. The extension class is represented by the polynomial
$p(z,u) = (p_{10} + p_{11}z + p_{12}z^2)u + p_{22}z^2u^2$. The generic
set of $\sM_3(2)$ is a projective $2$-space given by $[p_{10} : p_{11}
: p_{12}]$ minus the subvariety ${p_{10} = p_{12} = 0}$. Generic
bundles have invariants (width, height)$=(0,2)$, while the non-generic
bundles given by $u$ and $z^2u$ have invariants $(1,2)$, and the
bundle given by $z^2u^2$ has invariants $(2,2)$, \emph{like the split
bundle represented by $p=0$}. Thus the invariants do not separate
$\sM_3(2)$ into Hausdorff components (the split bundle is never
separated from any other bundle), which contrasts the situation for
instanton moduli $\sM_{nk}(k)$.
\end{exa}

\vfill

\noindent
Edoardo Ballico\\
University of Trento, Department of Mathematics\\
I--38050 Povo (Trento), Italia\\
E-mail: \url{ballico@science.unitn.it}

\bigskip
\bigskip
\noindent
Elizabeth Gasparim and Thomas K\"{o}ppe\\
School of Mathematics,
The University of Edinburgh\\
James Clerk Maxwell Building,
The King's Buildings,
Mayfield Road\\
Edinburgh, EH9 3JZ,
United Kingdom\\
E-mail: \url{Elizabeth.Gasparim@ed.ac.uk}\\
E-mail: \url{t.koeppe@ed.ac.uk}

\end{document}